\documentclass[11pt,a4paper,oneside,notitlepage,reqno]{amsart}\overfullrule=5pt

\usepackage{amsfonts}
\usepackage{amssymb}
\usepackage{amsthm}

\usepackage{verbatim} 
\usepackage{fancyhdr}
\numberwithin{equation}{section}

\DeclareMathSymbol{E}{\mathalpha}{AMSb}{"45}
\DeclareMathSymbol{\R}{\mathalpha}{AMSb}{"52}

\newcommand{\one}{\mathbf{1}}

\newcommand{\N}{\mathbb{N}}

\newcommand{\F}{\mathcal{F}}

\newcommand{\A}{\mathfrak{A}}

\newcommand{\ST}{\mathfrak{T}}

\renewcommand{\L}{\operatorname{L}}
\newcommand{\La}{\L^\alpha}
\newcommand{\Lah}{\L_h^\alpha}

\DeclareMathOperator{\D}{D}
\newcommand{\Dl}{\D_{\ell}}

\newcommand{\Dlk}{\D_{\ell_k}}

\DeclareMathOperator{\fddelta}{\delta}
\newcommand{\fddeltatau}{\fddelta_\tau}
\newcommand{\fddeltak}{\fddelta_{h,\ell_k}}
\newcommand{\fddeltamk}{\fddelta_{h,-\ell_k}}

\DeclareMathOperator{\fdDelta}{\Delta}
\newcommand{\fdDeltak}{\fdDelta_{h,\ell_k}}

\DeclareMathOperator{\T}{T}
\newcommand{\Thrlr}{\T_{h_r,\ell_r}}

\newcommand{\alphao}{{\alpha_0}}
\newcommand{\vrm}{v_r^-}

\newcommand{\mt}{\mathcal{M}_T}
\newcommand{\bmt}{\bar{\mt}}

\newcommand{\taut}{\tau_T(t)}

\newcommand{\on}{\,\, \textrm{on} \,\,}

\newcommand{\where}{\,\, \text{where} \,\,}
\newcommand{\cfor}{\,\, \text{for} \,\,}

\newcommand{\e}{\varepsilon}

\theoremstyle{plain}
\newtheorem{anytheorem}{Theorem}[section] 
\newtheorem{theorem}[anytheorem]{Theorem}
\newtheorem{lemma}[anytheorem]{Lemma}
\newtheorem{corollary}[anytheorem]{Corollary}
\theoremstyle{definition}
\newtheorem{anyassumption}{Assumption}[section]
\newtheorem{assumptions}[anyassumption]{Assumption}
\newtheorem{remark}[anytheorem]{Remark}
\newtheorem{example}[anytheorem]{Example} 

\title[On the rate of convergence]{On finite-difference approximations for normalized Bellman equations}
\author[I. Gy\"ongy]{Istv\'an Gy\"ongy}
\address{School of Mathematics and Maxwell Institute,
University of Edinburgh,
King's  Buildings,
Edinburgh, EH9 3JZ, United Kingdom}
\email{gyongy@maths.ed.ac.uk}
\author[D. Siska]{David  \v{S}i\v{s}ka}
\address{FIRST FRG, 
BNP Paribas,
10 Harewood Avenue,
London, NW1 6AA,
United Kingdom}
\email{davsiska@gmail.com}
\date{}

\begin{document}

\begin{abstract}
A class of stochastic optimal control problems  
involving optimal stopping is considered. 
Methods of 
Krylov \cite{krylov:rate:lipschitz:published} are
adapted to investigate the numerical solutions
of the corresponding normalized Bellman equations 
and to estimate the rate of convergence of  
finite difference approximations for the optimal 
reward functions. 
\end{abstract}


\maketitle

\section{Introduction}

Stochastic optimal control and optimal
stopping problems have many applications
in mathematical finance, portfolio
optimization, economics and statistics
(sequential analysis). Optimal stopping
problems can be in some cases solved
analytically
\cite{shiryaev:statistical}. With most
problems, one must resort to numerical
approximations of the solutions. One
approach is to use controlled Markov
chains as approximations to 
controlled diffusion processes, see
e.g. \cite{menaldi:some:estimates}. A
thorough account of this approach is
available in
\cite{kushner:dupuis:numerical}.

We are interested in the rate of
convergence of finite difference
approximations to the payoff function of
optimal stopping and control
problems. Using the method of randomized
stopping (see \cite{krylov:controlled}) 
such problems can be treated as optimal
control problems with the reward and
discounting functions unbounded in the
control parameter. This leads us to
approximating a normalized degenerate
Bellman equation.

Until quite recently, there were no
results on the rate of convergence of
finite difference schemes for degenerate
Bellman equations.  A major breakthrough
is achieved by Krylov in
\cite{krylov:rate:equations} for Bellman
equations with constant coefficients,
followed by rate of convergence
estimates for Bellman equations with
variable coefficients in
\cite{krylov:approximating:value} and
\cite{krylov:rate:variable}.  The
estimate from
\cite{krylov:rate:variable} is improved
in \cite{barles:jakobsen:error:bounds}
and
\cite{barles:jakobsen:rate:hamilton}.
Finally, Krylov
\cite{krylov:rate:lipschitz} (published
in
\cite{krylov:rate:lipschitz:published})
establishes the rate of convergence
$\tau^{1/4}+ h^{1/2}$ of finite
difference schemes to degenerate Bellman
equations with Lipschitz coefficients 
given on the whole space,
where $\tau$ and $h$ are the mesh sizes
in time and space respectively. This is 
later extended to 
finite difference approximations of 
Bellman equations on cylindrical
domains in \cite{dong:krylov:rate:domains}.

In the present paper we extend this
estimate to cover normalized degenerate
Bellman equations corresponding to
optimal stopping of controlled diffusion
processes with variable coefficients.
Adapting ideas and techniques of
\cite{krylov:rate:lipschitz} we obtain
the rate of convergence $\tau^{1/4} +
h^{1/2}$, as in
\cite{krylov:rate:lipschitz}.  
There are two 
key ideas which are already introduced 
in \cite{krylov:rate:equations}
-\cite{krylov:rate:variable}. The first idea is 
that the original equation and its approximation 
should play symmetric roles. The other idea is 
to `shake' the original equation and its 
approximation, and to mollify the solutions 
of the `shaken equations' to obtain smooth 
supersolutions to the original equation 
and to its approximation, respectively, 
which are close to their true solutions. 
To implement these ideas one needs 
appropriate estimates on the regularity 
of the solutions to the original equation 
and to its approximation. 
The necessary regularity estimates on the optimal 
reward functions, i.e., the solutions of 
the Bellman equations are well-known, see 
\cite{krylov:controlled}. Namely, under 
general conditions the optimal reward funtions 
are Lipschitz continuous in the space variable 
and they are H\"older continuous, with exponent $1/2$, 
in the time variable. The
main  problem 
is to obtain the corresponding regularity estimates 
for the finite difference approximations. 
In \cite{krylov:rate:lipschitz:published} 
a discrete gradient estimate in the space variable is
proved for the solutions to finite difference schemes  
for degenerate Bellman equations. Hence not only the 
Lipschitz continuity  in the space variable 
of the finite difference approximations follows but a
suitable estimate on their time regularity as well.

Our first main task in the present paper is to extend the 
discrete gradient estimate from
\cite{krylov:rate:lipschitz:published} 
to the case of finite difference schemes for 
normalized Bellman equations. 
This is Theorem \ref{thm-est-disc-grad} below. 
We note that in \cite{krylov:apriori} 
a more general estimate is proved. 
From Theorem \ref{thm-est-disc-grad} the Lipscitz 
continuity in the space variable of the finite difference 
approximations follows easily. 
However, due to the normalizing factor in the finite difference 
scheme, 
Theorem \ref{thm-est-disc-grad} does not imply 
the estimate we need on the time 
regularity  of the finite difference 
approximations. In fact, the time regularity 
of the solutions does not hold in general, unless we 
assume stronger conditions on the finite difference 
scheme than those of Theorem \ref{thm-est-disc-grad}. 
Since our main concern in the present paper 
is the rate of convergence of finite difference 
approximations for the reward function 
of optimal stopping of controlled diffusion processes, 
we establish the necessary time regularity estimate only for 
these approximations. This is Theorem \ref{cor:hold:cont:in:t}, 
which is the discrete counterpart of 
Theorem \ref{theorem 1.26.2.9} on the H\"older continuity 
in time of the optimal reward function. 
Hence, using also the regularity 
of the optimal reward functions 
and the maximal principle for normalized 
Bellman equations  
and for their `monotone approximations',
we prove our rate of convergence estimate, Theorem 
\ref{ass:opt:stop} by a straitforward adaptation of the  
method of `shaking and smoothing' 
from \cite{krylov:rate:lipschitz:published}.

Rate of converge results for optimal
stopping are proved for general
consistent approximation schemes in
\cite{jakobsen:rate:optimal:stopping}.
However, the rate $\tau^{1/4}+ h^{1/2}$
is obtained only when the diffusion
coefficients are independent of the time
and space variables.  For further
results on numerical approximations for
Bellman equations we refer to
\cite{jakobsen:karlsen:convergence:source:terms},
\cite{jakobsen:karlsen:chioma:error} and
\cite{biswas:jakobsen:karlsen:error}.

The paper is organized as follows.  The
main result, Theorem 
\ref{ass:opt:stop} is formulated in the next
section. In Section
\ref{section-existence-of-soln-to-disc-prob}
the existence and uniqueness of 
solutions to finite difference schemes, 
Theorem \ref{theorem 1.20.3.9}, is
proved together with a result 
on comparison of the solutions, 
Lemma \ref{lemma-discrete-comparison-principle}.
The
gradient estimate  
on the solutions of finite
difference schemes is proved in Section
\ref{section-gradient-estimate}, 
together 
with important corollaries.  
An estimate on 
Lipschitz continuity 
in the space variable for the reward functions  
and a result on comparison of the reward functions 
with supersolutions to Bellman equations are 
presented in
Section \ref{section-analytic-ppties}. 
The estimate on 
H\"older continuity in time  
of the reward  
functions together with the corresponding estimates 
for their finite difference approximations are 
given in Section \ref{section-hc}.
Theorem \ref{ass:opt:stop} is proved in
Section \ref{section-shaking}.

\section{The Main Result}                        \label{section-main-result}

Fix $T\in(0,\infty)$, 
and set $H_T=[0,T)\times\mathbb R^d$ 
and $\bar H_T=[0,T]\times\mathbb R^d$. 
Let $(\Omega, \F,P)$ 
be a probability space, carrying a 
$d'$ dimensional Wiener
martingale $W=(W_t)_{t\geq0}$ with respect 
to a filtration $(\F_t)_{t\geq0}$. 
Below we introduce some  
basic notions and notation of the 
theory of controlled diffusion 
processes from \cite{krylov:controlled}.  
The notation 
$|a|=(\sum_{i,j}a^2_{ij})^{1/2}$,  
$|b|=(\sum_{i}b_i)^{1/2}$ 
and $c^{+}=c_{+}=(|c|+c)/2$, $c^{-}=c_{-}=(-c)_{+}$ 
is used for 
matrices $a\in\mathbb R^{k\times l}$, vectors 
$b\in\mathbb R^k$ and real numbers $c$. Unless otherwise stated,  
the summation convention with respect 
to repeated indices is in force throughout the paper.

Let $A$ be a separable
metric space and let 
$\sigma=\sigma^{\alpha}(t,x)$, 
and $\beta=\beta^{\alpha}(t,x)$ be 
given Borel functions of 
$(\alpha,t,x)\in A\times\mathbb R\times\mathbb R^d$, 
taking values in $\mathbb R^{d\times d^{\prime}}$
 and  
$\mathbb R^d$, respectively. 
Assume that $A=\cup_{n=1}^{\infty}A_n$ for an 
increasing sequence 
of Borel sets $A_n$ of $A$ such that the following 
assumption holds.  

\begin{assumptions}                        \label{assumption 1.15.01.09}
For every integer $n\geq1$ there is a constant 
$K_n$ such that for all $\alpha\in A_n$ 
\begin{equation}                   \label{1.15.01.09}
|\sigma^\alpha(t,x) 
- \sigma^\alpha(t,y)| 
+ |\beta^\alpha(t,x) - \beta^\alpha(t,y)| 
\leq
K_n|x-y|
\end{equation}
\begin{equation}                   \label{2.15.01.09}
|\sigma^\alpha(t,x)| + |\beta^\alpha(t,x)|
\leq
K_n(1+|x|)
\end{equation}
for all $(t,x) \in \bar H_T$.
\end{assumptions} 

A progressively measurable process 
$\alpha = (\alpha_t)_{t\geq 0}$ 
with values in $A$ is called 
an {\it (admissible) strategy} if there is 
an integer $n\geq1$ such that 
$\alpha_t(\omega) \in A_n$ for all $t\geq0$ and $\omega\in\Omega$. 
The set of strategies with values in $A_n$ is denoted 
by $\mathfrak A_n$, and so 
$\mathfrak{A} =\bigcup_{n=1}^{\infty} \mathfrak{A}_n$ 
is the set of all strategies. 
By the classical existence and uniqueness theorem of 
It\^o, 
Assumption \ref{assumption 1.15.01.09} 
ensures that 
for each 
$\alpha \in \A$, $s\in[0,T]$ 
and  $x\in \R^d$ 
there is a unique solution
$x^{\alpha,s,x}=\{x_t:t\in[0,T-s]\}$ of
\begin{equation}                                         \label{eqn-sde}
x_t = x + \int_0^t \sigma^{\alpha_u}(s+u,x_u)dW_u 
+ \int_0^t \beta^{\alpha_u}(s+u,x_u)du.  
\end{equation}

Let $f=f^{\alpha}(t,x)$ and $c=c^{\alpha}(t,x)$ 
be Borel functions of 
$(\alpha,t,x)\in A\times\mathbb R\times\mathbb R^d$ 
with values in $\mathbb R$ and 
$\mathbb R_+$, respectively, and 
let $g=g(t,x)$ be a Borel function 
of $(t,x)\in\mathbb R\times\mathbb R^d$ with values in 
$\mathbb R$ 
 such that the following 
assumption holds. 

\begin{assumptions}                     \label{assumption 2.15.01.09}
The function $g$ is continuous 
and there are some constants
$K$ and
$q\geq0$  such that 
\begin{equation}                         \label{7.15.01.09}
|g(t,x)| \leq K(1+|x|^q) 
\quad \text{for all $(t,x)\in\bar H_T$.}
\end{equation}
For every integer 
$n\geq1$ there  are  constants 
$K_n$ and $q_n\geq0$ such that 
for all $\alpha\in A_n$
\begin{equation}                          \label{8.15.01.09}
c^\alpha(t,x)\leq K_n(1+|x|^{q_n}), 
\quad  |f^\alpha(t,x)| 
\leq K_n(1+|x|^{q_n})
\end{equation}
for all $(t,x)\in \bar H_T$.
\end{assumptions}

 For
$s\in[0,T]$ we use the notation
$\ST(T-s)$ for the set of stopping times
$\tau\leq T-s$.  Consider the  
following {\it optimal reward functions}: 
\begin{equation}                                  \label{3.17.01.09}
v(s,x)=\sup_{\alpha\in\mathfrak{A}}v^{\alpha}, 
\quad (s,x)\in \bar H_T, 
\end{equation}
\begin{equation}                                 \label{eq:stop-payofffn}
w(s,x)=\sup_{\alpha\in\mathfrak{A}}
\sup_{\tau\in\ST(T-s)}
w^{\alpha,\tau}(s,x), 
\quad (s,x)\in \bar H_T, 
\end{equation}
where 
\begin{equation}                               \label{5.17.01.09}
v^{\alpha}(s,x)=E^\alpha_{s,x}
\left[\int_0^{T-s}f^{\alpha_t}(s+t,x_t)
e^{-\varphi_t}dt + g(T,x_{T-s})e^{-\varphi_{T-s}}
\right],
\end{equation}
\begin{equation}                               \label{1.03.02.09}
w^{\alpha,\tau}(s,x)=E^\alpha_{s,x}
\left[\int_0^{\tau}f^{\alpha_t}(s+t,x_t)
e^{-\varphi_t}dt + g(s+\tau,x_{\tau})e^{-\varphi_{\tau}}
\right],
\end{equation}
$$                                         
\varphi_t = \varphi_t^{\alpha,s,x}
= \int_0^t c^{\alpha_r}(s+r,x_r^{\alpha,s,x})dr, 
$$
and $E^\alpha_{s,x}$ denotes the expectation
of the expression behind it, with
$x^{\alpha,s,x}_t$ in place of $x_t$
everywhere.  
We call $v$ and $w$  the optimal reward functions 
for the {\it optimal  control} problem, and for the 
{\it optimal control and stopping} problem, respectively, 
with strategies from $\A$, under 
{\it utility rate} $f$, {\it terminal utility} $g$ 
and {\it discount rate} $c$.
It is useful to notice that for
$$
v_n(s,x)=\sup_{\alpha\in\mathfrak{A}_n}v^{\alpha}(s,x), 
\quad 
w_n(s,x):=\sup_{\alpha\in\mathfrak{A}_n}
\sup_{\tau\in\ST(T-s)}
w^{\alpha,\tau}(s,x)
$$
we have $v_n(s,x)\uparrow v(s,x)$ and $
w_n(s,x)\uparrow w(s,x)$ as 
$n\to \infty$.  
Our aim is to investigate 
finite difference approximations  
for a class of nonlinear PDEs, called 
{\it normalized Bellman PDEs}, to approximate 
$w$ via finite difference schemes for appropriate 
normalized Bellman PDEs,  
and to study the accuracy of these approximations.   

Using the method 
of {\it randomized stopping}, it is very useful 
to rewrite 
\eqref{eq:stop-payofffn} 
in the form of \eqref{3.17.01.09}, 
by extending $A_n$ and $\A_n$ as follows.  
Set 
$$ 
\bar A=A\times[0,\infty)=\cup_{n=1}^{\infty}\bar A_n, 
\quad 
\bar A_n=A_n\times[0,n],    
$$
identify $\alpha \in A$ with $(\alpha,0)\in\bar A$, 
and extend the definition of 
$\sigma$,
$\beta$,
$f$,
$g$ and
$c$  by setting  
$$
\sigma^{\gamma}=\sigma^{\alpha},
\quad
\beta^{\gamma}=\beta^{\alpha},
\quad
f^{\gamma}=f^{\alpha}+rg, 
\quad
c^{\gamma}=c^{\alpha}+r,
\quad\text{for $\gamma=(\alpha,r)\in\bar A$}.
$$ 
Let $\bar\A_n$ denote the set of progressively 
measurable processes with values in $\bar \A_n$ 
and set $\bar \A=\cup_n\bar\A_n$.  
Notice, that if Assumptions 
\eqref{assumption 1.15.01.09}-\eqref{assumption 2.15.01.09}
hold then these assumptions remain 
valid with $\bar A_n$ and $\bar A$ in place 
of $A_n$ and $A$, with the obvious extension 
of the metric on $A$ onto $\bar A$.  
Moreover, the following result holds.

\begin{theorem}                        \label{thm-opt-stop-and-cont-equiv}
  Let Assumptions
  \ref{assumption 1.15.01.09} 
and \ref{assumption 2.15.01.09} 
hold. Then $w=\sup_{\gamma\in\bar\A}v^{\gamma}$ 
for every $(s,x) \in [0,T]$, 
where $v^{\gamma}$ is defined by 
\eqref{5.17.01.09} with $\gamma\in\bar\A$ 
in place of $\alpha\in\A$. 
\end{theorem}
This theorem, under somewhat stronger assumption 
is known from 
\cite{krylov:controlled} 
when 
$A=A_n$, $K=K_n$, $m=m_n$ for $n\geq1$. 
For the proof we refer to 
\cite{gyongy:siska:on:randomized}.

From \cite{krylov:controlled} one also knows 
that under some assumptions (more
restrictive than Assumptions
\ref{assumption 1.15.01.09}-\ref{assumption 2.15.01.09})
$w$ satisfies
the {\it normalized Bellman PDE} 
\begin{equation}                            \label{eq:chp:rate:3}
    \sup_{\gamma\in\bar A}
m^{\gamma}( 
\tfrac{\partial}{\partial t}w 
+ L^{\gamma} w +
f^{\gamma})  = 0 \on H_T 
\end{equation}
with terminal condition
\begin{equation}                            \label{6.20.01.09}         
w(T,x) = g(T,x) \cfor x\in \R^d,
\end{equation} 
where 
$m^{\gamma}=(1+r)^{-1}$ 
and
\begin{equation}                                      \label{eqn-La}
  L^{\gamma} w = \tfrac{1}{2}\sigma^{\gamma}_{ip}\sigma^{\gamma}_{jp}
w_{x^i x^j} + 
  \beta^{\gamma}_i w_{x^i} - c^{\gamma} w.  
\end{equation} 
Therefore it is natural to design 
approximations for $w$ as  
finite difference approximations  
for problem \eqref{eq:chp:rate:3}-\eqref{6.20.01.09}.  
To this end we fix a constant 
$K\geq1$ and 
make the assumptions below. 

\begin{assumptions}                         \label{ass-on-the-scheme}
There exist a natural number $d_1$,
vectors $\ell_k \in \R^d$ and functions
$$
a^\alpha_k:\mathbb R\times \R^d \to \R_{+}, 
\quad b^\alpha_k:\mathbb R\times \R^d \to \R_{+}, 
\quad \text{for $k = \pm 1, \ldots,\pm d_1$ and 
$\alpha\in A$},
$$
such that   
$|\ell_k| \leq K$, 
$\ell_k = -\ell_{-k}$, 
$a^{\alpha}_k =a^{\alpha}_{-k}$, 
  for $k = \pm 1, \ldots, \pm d_1$, 
$\alpha \in A$, and 
\begin{equation}                     \label{1.01.19.09}
    \beta^\alpha_i =  b_k^\alpha \ell_k^i 
\end{equation}
\begin{equation}                   \label{2.01.19.09}
    \tfrac{1}{2}\sigma^{\alpha}_{ip} \sigma^{\alpha}_{jp}
    = a_k^\alpha \ell_k^i \ell_k^j, 
\end{equation}
for $\alpha\in A$ and  $i,j=1,2,\dots,d$. 
\end{assumptions}

\begin{remark}
For given functions $\beta^{\alpha}$ it is easy to 
find a set of vectors $\{\ell_k\}$ and functions  $b^{\alpha}_k\geq0$ 
such that \eqref{1.01.19.09} holds. We can take,  for example,  
$\ell_{\pm k}=\pm e_k$, with the standard basis $\{e^k\}$ in 
$\mathbb R^d$, and set $b_{\pm k}^{\alpha}=(\beta_k^{\alpha})_{\pm}$. 
It is proved in \cite{krylov:factorizations} that, if 
the matrix $\sigma^{\alpha}\sigma^{\alpha\ast}$
is uniformly nondegenerate, then there always exist 
a set of vectors 
$\ell_k\in\mathbb R^d\setminus\{0\}$ and 
functions $a_k^{\alpha}$ for $k=\pm1,\dots,\pm d_1$ 
for some integer $d_1$ such that 
$\ell_{-k}=-\ell_k$, $a^{\alpha}_{-k}=a^{\alpha}_k\geq0$ 
for all $k$, \eqref{2.01.19.09} holds,
 $a^{\alpha}_{k}$ are as smooth as 
$\sigma^{\alpha}\sigma^{\alpha\ast}$ is,
and $a_{k}^{\alpha}\geq\kappa>0$,   
where $\kappa$ is a constant.
It is also proved  in \cite{krylov:factorizations} that if 
all values of the matrix $\sigma^{\alpha}\sigma^{\alpha\ast}$ 
lie in a closed convex
polyhedron in the set of nonnegative matrices and  
the first and second order derivatives in $x\in\mathbb R^d$
of $\sigma^{\alpha}\sigma^{\alpha\ast}$ 
are bounded  functions, then again
there   exist
$\{\ell_k\}$ and $a_k^{\alpha}$
satisfying the above assumption such that 
$\sqrt{a_k^{\alpha}}$
are Lipschitz continuous in $x$. 

\end{remark}
Clearly, \eqref{1.01.19.09} and \eqref{2.01.19.09} imply 
$$
\tfrac{1}{2}\sigma^{\alpha}_{ip}\sigma^{\alpha}_{jp}u_{x^ix^j}
=a_k^{\alpha}D_{\ell_k}^2u, 
\quad 
\beta_i^{\alpha}u_{x^i}=b_k^{\alpha}D_{\ell_k}u
$$
for smooth functions $u$, where we use the notation 
$$
\quad \Dl u = u_{x^i}\ell^i 
\quad
\text{for $\ell\in\mathbb R^d$}. 
$$
Thus setting 
$a^{\gamma}_k=a^{\alpha}_k$ 
and $b^{\gamma}_k=b^{\alpha}_k$ for 
$\gamma=(\alpha,r)\in\bar A$, 
for the operator $L^{\gamma}$ given 
by \eqref{eqn-La} we have 
\begin{equation*}                           \label{eqn-Lak}
   L^{\gamma} u =  a^{\gamma}_k 
\Dlk^2 u + b^{\gamma}_k \Dlk u 
- c^{\gamma}u,
\quad\text{for $\gamma\in\bar A$}. 
\end{equation*}

For $\tau > 0$, $h > 0$ and $l\in\mathbb R^d$
define  
\begin{equation}                                        \label{eqn-f-d-ops}
  \begin{split}
    \fddeltatau u(t,x) 
& := \tfrac{u(t+\tau,x) -
u(t,x)}{\tau},  
\quad
\tau_T(t)=\tau\wedge(T-t)                       \\
    \fddeltatau^T u(t,x) 
&:= \tfrac{u(t+\tau_T(t),x) -
u(t,x)}{\tau},                                   \\
    \delta_{h,l} u(t,x) 
& := \tfrac{u(t,x+hl) -
u(t,x)}{h},                                       \\
    \Delta_{h,l} u &
:= - \delta_{h,l}\delta_{h,-l} u 
= \tfrac{1}{h}(\delta_{h,l}u + \delta_{h,-l}u).
  \end{split}
\end{equation}
for $t\in[0,T)$, $x\in\mathbb R^d$, and consider  
the finite difference scheme

\begin{equation}                           \label{1.20.01.09}
\sup_{\gamma\in\bar A}m^{\gamma} 
(\fddeltatau^{T} u +L^{\gamma}_hu 
+ f^{\gamma}) =  0 
\quad \text{on $H_T$}
\end{equation}  
\begin{equation}                          \label{7.20.01.09}
    u(T,x)  =  g(x)
\quad \text{for $x\in\mathbb R^d$},
\end{equation}
where
\begin{equation*}
  L_h^{\gamma}u = a_k^{\gamma} \fdDeltak u + b_k^{\gamma}
\fddeltak u - c^{\gamma} u.
\end{equation*}

\begin{remark}                              \label{remark 7.20.01.09}
Equation \eqref{eq:chp:rate:3} is often written in 
the form 
\begin{equation}                           \label{1.01.02.09}
\max\big[
\tfrac{\partial}{\partial t} w
 +\sup_{\alpha\in A}(L^{\alpha}w 
+ f^{\alpha}) ,g-w\big]=0\quad \text{on $H_T$}, 
\end{equation}
and similarly, equation \eqref{1.20.01.09} can be written as
\begin{equation}                            \label{2.01.02.09}
\max\big[\fddeltatau^{T} u  
+\sup_{\alpha\in A}(L^{\alpha}_h u
+ f^{\alpha}) ,g-u \big]=0\quad \text{on $H_T$}. 
\end{equation}
Clearly, equation \eqref{1.01.02.09} 
is equivalent to 
$$
\tfrac{\partial}{\partial t}w 
+ \sup_{\alpha \in A}\left[L^{\alpha}w + f^\alpha \right] 
\leq  0,  \quad g - w \leq 0 
\quad \text{on $H_T$},
$$
$$
\tfrac{\partial}{\partial t}w 
+ \sup_{\alpha \in A}\left[L^{\alpha}w + f^\alpha \right] 
= 0,  \quad \text{on $\{(t,x)\in H_T:g(t,x)< w(t,x)\}$}, 
$$
and similarly 
equation \eqref{2.01.02.09} is equivalent to 
$$
\delta_{\tau}^{T}u
+ \sup_{\alpha \in A}\left[\Lah u + f^\alpha \right]  
\leq  0, \quad g - u \leq 0 \quad \text{on $H_T$},
$$
$$ 
\delta_{\tau}^{T}u  + \sup_{\alpha \in A} 
\left[ \Lah u  + f^\alpha \right] = 0 \quad 
\text{on $\{(t,x)\in H_T:g(t,x)< u (t,x)\}$}.   
$$
\end{remark}

\begin{proof}
By setting $\e = \tfrac{1}{1+r}$ equations 
\eqref{eq:chp:rate:3} and 
\eqref{1.20.01.09} can be rewritten as 
$$
 \sup_{\e \in [0,1]}
\left[ \e \sup_{\alpha \in A} 
\left(\tfrac{\partial}{\partial t}w  
+L^{\alpha}w + f^\alpha \right) + (1-\e)(g-w) \right] = 0
\quad\text{on $H_T$}                                                 
$$
and 
$$
\sup_{\e \in [0,1]}
\left[ \e \sup_{\alpha \in A} 
\left( \delta_{\tau}^{T} u
+L_{h}^{\alpha} u  + f^\alpha \right) + (1-\e)(g-u) \right] =0
\quad\text{on $H_T$}, 
$$
respectively. 
Hence we finish the proof of the remark by noticing 
that for any numbers $p,q\in\mathbb R$ 
$$
\sup_{\e \in [0,1]}(\varepsilon p+(1-\varepsilon) q)=\max(p,q).  
$$
\end{proof}

Instead of Assumptions \ref{assumption 1.15.01.09} 
and \ref{assumption 2.15.01.09} we make now the following 
assumption. 

\begin{assumptions}                        \label{ass:opt:stop}
 The functions 
$\sigma^{\alpha}$, $a^{\alpha}_k$, 
  $b^\alpha_k$, $f^\alpha$ and $c^\alpha\geq 0$ 
are Borel measurable in $t$ and 
are continuous in $\alpha \in A$ for
  each $k = \pm 1, \ldots, d_1$. Moreover, 
for $\Psi:=\sigma^{\alpha},\sqrt{a^{\alpha}_k}, b_k^{\alpha}, 
c^{\alpha}, f^{\alpha}, g$ for $\alpha\in A$ and 
$k=\pm 1, \ldots, \pm d_1$ we have 
  \begin{equation}                          \label{eq:bds:opt:stop}
 |\Psi(t,x)-\Psi(t,y)|\leq K|x-y|,
\quad |\psi(t,x)|\leq K
  \end{equation}
for all $t\in\mathbb R$ and $x\in\mathbb R^d$. 
\end{assumptions}
Notice that Assumption 
\ref{ass-on-the-scheme} and 
\ref{ass:opt:stop} imply 
Assumptions \ref{assumption 1.15.01.09} 
and \ref{assumption 2.15.01.09}. 
Finally we make an assumptions on H\"older
continuity of $\sqrt{a^{\alpha}_k}$,
$b^\alpha_k$, $c^\alpha$ and $f^\alpha$.

\begin{assumptions}                     \label{ass-for-main-result}
For $\Psi:=\sqrt{a^{\alpha}_k},
b^\alpha_k, c^\alpha,f^\alpha,g$ 
for $\alpha\in A$, 
$k=\pm 1, \ldots, \pm d_1$ we have  
  \begin{equation*}
 |\psi(t,x)-\psi(t,x)|\leq K|t-s|^{1/2}
  \end{equation*}
for all $x\in\mathbb R^d$ and $s,t\in\mathbb R$. 
\end{assumptions}

The following result is the main theorem of the paper. 
It extends 
Theorem 2.3 
from \cite{krylov:rate:lipschitz:published}  
to the reward function $w$ 
defined by \eqref{eq:stop-payofffn}. 

\begin{theorem}                          \label{theorem-the-main-result}
Let Assumptions 
\ref{ass-on-the-scheme} through 
\ref{ass-for-main-result} hold.  
Then 
\eqref{1.20.01.09}-\eqref{7.20.01.09}
has a unique bounded
solution $w_{\tau,h}$, and there is a 
constant $N$ depending only on $K,d,d_1,T$ such that 
for $\tau,h\leq 1$ 
\begin{equation}                \label{5.20.01.09}
    |w - w_{\tau,h}| 
\leq N(\tau^{1/4} + h^{1/2}) 
\end{equation}
on $\bar H_T$. Moreover, there is a constant 
$\lambda_0$ depending only on $K$ and $d_1$ such that if 
$\lambda\geq\lambda_0$ then $N$ is independent of $T$. 

\end{theorem}

\section{On finite difference schemes}
                               \label{section-existence-of-soln-to-disc-prob}

Let $A$ be a set and consider for $\alpha\in A$ 
the finite difference
operator 
\begin{equation*}
  \Lah = a_k^\alpha \fdDeltak  + b_k^\alpha
\fddeltak  - c^\alpha ,
\end{equation*}
where $a^\alpha_k$, $b^\alpha_k$, $c^\alpha$, $f^\alpha$ 
and $g$ are some functions on $H_{\infty}:=[0,\infty)\times\R^d$ 
for each $\alpha \in A$ and $k = \pm 1,\ldots, \pm d_1$. 
Recall that $\{\ell_k:k=\pm1,\pm2,\dots,\pm d_1\}$ 
are given vectors in  $\mathbb R^d$ such that
$|\ell_k|
\leq K$ for all $k = \pm 1, \ldots, \pm d_1$ and 
$l_k=-l_k$, where $K\geq 1$ is a fixed constant.  

Let $m^{\alpha}$ be a function of $\alpha\in A$ 
taking values in $(0,1]$. 
Recall that $H_T=[0,T)\times\mathbb R^d$ 
for a fixed $T\in[0,\infty)$.  
For fixed $\tau>0$ and $h>0$ we are interested 
in the problem 
\begin{equation}                                               \label{2.30.12.08}
  \sup_{\alpha \in A} 
m^\alpha \left(\fddeltatau^T v 
+ \Lah v + f^\alpha \right)  
=  0 \quad \textrm{on}\quad H_T,
\end{equation}
\begin{equation}                                                 \label{3.30.12.08}
  v(T,x) =  g(T,x) \quad x\in\mathbb R^d                   
\end{equation}
for a function $v=v_{\tau,h}$ defined on 
$\bar H_T=[0,T]\times\mathbb R^d$. 
Notice that problem 
\eqref{2.30.12.08}-\eqref{3.30.12.08} 
is a collection of separate problems 
given on each grid
\begin{equation}                                                   \label{4.30.12.08}
\{((t_0+j\tau)\wedge T,x_0+(\pm i_1\ell_1\pm\dots\pm i_{d_1}\ell_{d_1})h\}
\end{equation}
associated with points 
$(t_0,x_0)\in[0,T)\times\mathbb R^d$, 
where $i_1$,...,$i_{d_1}$ and $j$ 
run through the nonnegative integers. 
The grid associated with 
the point
$(t_0,x_0):=(0,0)$ is 
$$
\bar{\mathcal{M_T}}=\{(j\tau\wedge T,\pm i_1h\ell_1\pm\dots\pm
i_{d_1}h\ell_{d_1}):j,i_1,\dots,i_{d_1}=0,1,\dots\}. 
$$ 
Clearly, results obtained for equations on subsets of 
$$
\mathcal M_T:=\bar{\mathcal M_T}\cap([0,T)\times\mathbb R^d)
$$ 
can be translated into
results  for for equations on subsets  of all other grids 
of the type \eqref{4.30.12.08}.

In this section we consider the 
finite difference problems
\begin{equation}                                              \label{eq:fd:a}
  \sup_{\alpha \in A} 
m^\alpha \left(\fddeltatau^T u 
+ \Lah u + f^\alpha \right)  
=  0 \quad \textrm{on}\quad Q,
\end{equation}
\begin{equation}                                               \label{eq:fd:b}
  u =  g \on \bar{\mathcal{M}}_T\setminus Q                    
\end{equation}
and 
\begin{equation}                                              \label{1.17.02.09}
 \max\left[\sup_{\alpha \in A} 
m^\alpha (\fddeltatau^T w 
+ \Lah w + f^\alpha ),g-w\right]  
=  0 \quad \textrm{on}\quad Q,
\end{equation}
\begin{equation}                                               \label{2.17.02.09}
  w =  g \on \bar{\mathcal{M}}_T\setminus Q                    
\end{equation}
where $Q$ is a fixed subset of 
$\mathcal M_T$ and 
$g$ is a bounded function on 
$H_{\infty}$. 
Let $\lambda\geq0$ be a constant and make the 
following assumptions. 

\begin{assumptions}                                   \label{assumption 1.23.12.08}
We have $m^{\alpha}\in(0,1]$, 
$a^\alpha_k \geq 0$, $b^{\alpha}_k\geq0$, 
$a^\alpha_k = a^\alpha_{-k}$
and $c^\alpha \geq \lambda$ for all 
$\alpha\in A$, $(t,x)\in H_{\infty}$ and 
$k=\pm1,\pm2,\dots, \pm d_1$. 
\end{assumptions}
\begin{assumptions}                                       \label{ass-new-for-3}
For all $k = \pm 1, \ldots, \pm d_1$, 
$\alpha \in A$, $(t,x) \in H_{\infty}$
\begin{equation*}                                        \label{ass-new-for-3.1}
|m^\alpha a^\alpha_k| 
+ |m^\alpha b^\alpha_k| + |m^\alpha c^\alpha| 
+ |m^\alpha f^\alpha|\leq K. 
\end{equation*}
\end{assumptions}
\begin{assumptions}                                   \label{assumption 2.23.12.08}
There exists a constant $\rho>0$ 
such that 
\begin{equation}                                           \label{eqn-ass-m-c}
  m^\alpha(1+c^\alpha-\lambda) \geq \rho
\end{equation} 
on $H_{\infty}$ for all $\alpha\in A$.                               
\end{assumptions}

\begin{remark}                                          \label{remark 3.17.02.09}
Consider $\bar A=A\times[0,\infty)$, 
identify every $\alpha\in A$ with $(\alpha,0)\in \bar A$,  
and set for $\gamma=(\alpha,r)\in \bar A$
$$
m^{\gamma}=m^{\alpha}(1+r)^{-1}, \quad 
a_k^{\gamma}=a_k^{\alpha}, \quad
b_k^{\gamma}=b_k^{\alpha}, 
$$
$$
c^{\gamma}=c^{\alpha}+\tfrac{r}{m^{\alpha}}, 
\quad f^{\gamma}=f^{\alpha}+\tfrac{r}{m^{\alpha}}g. 
$$
$$
L^{\gamma}_h = a_k^{\gamma} \fdDeltak  + b_k^{\gamma}
\fddeltak  - c^{\gamma} . 
$$
Then, as Remark \ref{remark 7.20.01.09} is shown, 
it is easy to see that 
equation \eqref{1.17.02.09} can be cast into 
equation \eqref{eq:fd:a} with $\bar A$ in place of $A$. 
Clearly, if Assumption \ref{assumption 1.23.12.08} holds, 
then it holds also with $\bar A$ in place of $A$. 
If Assumption \ref{assumption 2.23.12.08} holds, 
then it is easy to show that 
it holds with $\bar A$ in place of 
$A$ and with $\min(\rho,1)$ in place of $\rho$. 
If Assumption \ref{ass-new-for-3} 
holds and $|g|\leq K$ on $H_{\infty}$ 
then it is easy to see that 
Assumption \ref{ass-new-for-3} holds also with $\bar A$ 
in place of $A$,  
with constant $2K+1$ in place of $K$. 
Thus we obtain the results of this 
section immediately for both equations 
\eqref{eq:fd:a} and \eqref{1.17.02.09},  
by proving them only for \eqref{eq:fd:a} 
and verifying that the conditions formulated 
with $A$ hold also with $\bar A$ in place of $A$.  
\end{remark}

The following simple 
examples show that if condition 
\eqref{eqn-ass-m-c} does not hold then 
problem \eqref{eq:fd:a}-\eqref{eq:fd:b} 
may have many solutions or may have no 
solution. 
\begin{example}
Let $A=[0,\infty)$, $m^\alpha =
  (1+\alpha)^{-1}$. Consider  the problem 
$$
    \sup_{\alpha \in A} m^\alpha
    \left(\fddeltatau^T u \right) = 0 \on
    \mt, 
\quad u = 1 \quad \text{on $\bmt \setminus \mt$}.  
$$
Notice that here 
  $\inf_{\alpha \in A} m^\alpha(1+c^\alpha) = 0$,  
  i.e. the condition 
  \eqref{eqn-ass-m-c} is violated.
  If $u:\mt \to \R$
  is any non-increasing function in $t$,
  then $m^\alpha \fddeltatau^T u \leq
  0$. Hence, letting $\alpha \to \infty$,
  we see that $u$ satisfies the
  equation. Consequently the solution to
  the above problem is not unique.
\end{example}
\begin{example}
  Let $A=[0,\infty)$, $m^\alpha =
  (1+\alpha)^{-1}$ and $f^\alpha = 1 + \alpha$. 
  Consider now the equation
  \begin{equation*}
    \sup_{\alpha \in A}m^\alpha(\fddeltatau^T u 
+ f^\alpha) = \sup_{\alpha \in A}m^\alpha\fddeltatau^T u + 1 
= 0 \on \mt.
  \end{equation*}
  If $u$ is a solution then we have
  $m^\alpha \fddeltatau^T u \leq 0$. Hence
  $\sup_{\alpha \in A}m^\alpha
  \fddeltatau^T u = 0$, which contradicts
  the equation. Thus the above equation
  has no solution.
\end{example}

\begin{theorem}               \label{theorem 1.20.3.9}
  Let Assumptions 
\ref{assumption 1.23.12.08} through 
\ref{assumption 2.23.12.08} hold. 
Let $g$ be a bounded function on 
$\bar{\mathcal M}_T$. Then  
  the finite difference problems 
  \eqref{eq:fd:a}-\eqref{eq:fd:b} 
and  \eqref{1.17.02.09}-\eqref{2.17.02.09} 
admit a unique bounded solution $u$ and $w$, respectively. 
\end{theorem}

\begin{proof}
By virtue of Remark \ref{remark 3.17.02.09} 
it suffices to prove the lemma for 
\eqref{eq:fd:a}-\eqref{eq:fd:b}. 
Let $\gamma = (0,1)$ and define $\xi$
  recursively as follows: $\xi(T) = 1$,
  $\xi(t) = \gamma^{-1} \xi(t + \taut)$
  for $t < T$. Then for any function $v$
  \begin{equation*}
    \fddeltatau^T
(\xi v) = \gamma \xi \fddeltatau^T v - \nu \xi
v, 
\where \nu = \tfrac{1-\gamma}{\tau}.
  \end{equation*}
 Solving 
  \eqref{eq:fd:a}-\eqref{eq:fd:b}
  for $u$ is equivalent to solving 
  \begin{equation}\label{eqn-proxy-for-eqn-f-d-problem}
    v = H[v] := H[(f^\alpha),g,v] :
    = \one_{\bar{\mathcal{M}}_T\setminus Q}\tfrac{1}{\xi}g + \one_{Q}G[v], 
  \end{equation}
  with $u = \xi v$, 
  where for $\e > 0$,
  \begin{equation}                                             \label{eq:G}
    G[v] := v + \e \xi^{-1} \sup_{\alpha}
    m^\alpha \left( \fddeltatau^T u 
      + \Lah u + f^\alpha \right).
  \end{equation}
  Then
  \begin{equation}                                            \label{eq:G2}
    G[v] = \sup_\alpha \left[ p^\alpha_\tau \T_\tau v 
      + p^\alpha_k \T_{h,l_k} v +
      p^\alpha v + \e m^\alpha 
      \xi^{-1}f^\alpha\right],
  \end{equation}
  with
  \begin{equation*}
    p^\alpha_\tau = \e \gamma \tau^{-1}
 m^\alpha \geq 0, 
\quad p^\alpha_k 
= \e(2 h^{-2}a^\alpha_k + h^{-1}b^\alpha_k)m^\alpha \geq 0,
  \end{equation*}
  \begin{equation*}
    p^\alpha = 1 - p^\alpha_\tau - \sum_k p^\alpha_k 
- \e \nu m^\alpha - \e m^\alpha c^\alpha.
  \end{equation*}
  Notice that 
  $
p^\alpha_k \leq \e K (h^{-2} + h^{-1})
$
  and
$$ 
\e \nu m^\alpha + \e m^\alpha c^\alpha 
\leq \e \tau^{-1} + \e K, 
\quad p_\tau^\alpha
 \leq \e \tau^{-1},
$$
  so for all $\e$ smaller than some
  $\e_0$ we have $p^\alpha \geq 0$. Also 
by taking into account \eqref{eqn-ass-m-c}
we have
\begin{equation*}
    \begin{split}
      0 \leq \sum_kp^\alpha_k + p^\alpha 
      + p_\tau^\alpha & 
      = 1 - \e m^\alpha (\nu + c^\alpha) 
      \leq 1 - \e (1 \wedge \nu)
      m^\alpha(1+ c^\alpha) \\
      & \leq 1 - \e (1 \wedge \nu)\rho
      =: \delta < 1,
    \end{split}
\end{equation*}
  for sufficiently small $\e >
  0$. Notice also $|m^\alpha f^\alpha| \leq K$. Hence
  $H$ maps bounded functions on $\bmt$ into bounded 
  functions on $\bmt$. Furthermore
\begin{equation*}
    |H[v](t,x) - H[w](t,x)| 
\leq \delta \sup_{\bar{\mathcal{M}}_T}|v-w|.
\end{equation*}
  Thus the operator $H$ is a contraction
  on the space of bounded functions on
  $\bar{\mathcal{M}}_T$. By Banach's fixed
  point theorem
  \eqref{eqn-proxy-for-eqn-f-d-problem}
  has a unique bounded solution.
\end{proof}
Set 
$\bar{\mathcal M}_{T,R}
=\{(t,x)\in\bar{\mathcal M_T}, |x|\leq R\}$ 
and 
$\bar{\mathcal M}_{T,R}^c
=\{(t,x)\in\bar{\mathcal M}_T, |x|>R\}$ 
for $R>0$.
\begin{remark}                                        \label{rem-frp}
  Let $v$ be a function on
  $\bar{\mathcal M_T}$. The operator $H$ defined by
  \eqref{eqn-proxy-for-eqn-f-d-problem}   
has the following property: if there   
exists $R > 0$ such that $v=f^\alpha = 0$ 
on $\bar{\mathcal M}_{T,R}^c$
for all   $\alpha \in A$, then there exists $R'$
such that   
\begin{equation*}
    H[(f^\alpha),g,v](t,x) = 0
\quad\text{on $\bar{\mathcal M}_{T,R'}^c$}.  
\end{equation*}
\end{remark}
\begin{corollary}                                       \label{corollary-R}
  Let Assumptions 
\ref{assumption 1.23.12.08} through 
\ref{assumption 2.23.12.08} 
hold. 
  Let $u$ be the bounded solution of
  \eqref{eq:fd:a}-\eqref{eq:fd:b}
  with $Q = \mt$. Assume there exists $R>0$ such
  that for all $\alpha \in A$ 
  $$
f^\alpha = g= 0 
\quad\text{on $\bar{\mathcal M}_{T,R}^c$}.
$$ 
  Then 
$$
\lim_{r\to \infty} \sup_{\bar{\mathcal M}_{T,r}^c} |u(t,x)| = 0.
$$
\end{corollary}
\begin{proof}
  Let $\xi$ be defined as in the proof
  of Theorem 
  \ref{theorem 1.20.3.9}
  and let $v = \xi u$. 
  For a fixed $(f^\alpha)$ and $g$ we define $H^n[v]$ 
  for functions $v$ on $\bmt$ 
  recursively in $n$ as follows: $H^1[v] = H^1[(f^\alpha), g, v]$
  and $H^n[v] = H^1[H^{n-1} [v]]$ for $n\geq2$.
  From the proof of
  Theorem
  \ref{theorem 1.20.3.9}
  we see that $H$ is a contraction on   
  the space of bounded functions on   $\bmt$. 
  Hence for any $\e > 0$ there   is $n_0$
  such that
  \begin{equation*}
    \sup_{\bar{\mathcal M}_T}|H^{n}[0] - v| < \e, \cfor n \geq n_0.
  \end{equation*}
  By Remark \ref{rem-frp} there exist
  $R_\e$ such that $H^{n_0}[0]=0$
  on $\bar{\mathcal M}_{T,R_{\varepsilon}}^c$. Hence  
  \begin{equation*} 
    \sup_{\bar{\mathcal M}_{T,R_{\varepsilon}}^c}|v| < \e,  
  \end{equation*}
which proves the corollary. 
\end{proof}

For the next lemma we need some   
remarks from \cite{krylov:rate:lipschitz}.
Let $\D_x^n$ denote the collection of all
n-th order derivatives in $x$. 

\begin{remark}                               \label{remark-taylors-thm}
For any sufficiently smooth function 
$\eta(x)$ by Taylor's formula
\begin{equation*}
|\La \eta(x) - \Lah \eta(x)| 
\leq N(h^2 \sup_{B_K(x)}|\D^4_x\eta| + h \sup_{B_K(x)} |\D^2_x \eta|), 
\end{equation*}
where $B_K(x)$ is the ball of radius $K$ centered 
at $x$.
\end{remark}

\begin{remark}                                \label{remark-T'}
  Let us introduce $T'$ as the least integer
  multiple of $\tau$ not less than $T$. 
  Notice that problem \eqref{eq:fd:a}-\eqref{eq:fd:b} can be rewritten
  as 
  \begin{equation*}
    \sup_{\alpha \in A}\left(\fddeltatau \tilde{u} 
      + \Lah \tilde{u} + f^\alpha\right) = 0,
\quad\text{on $Q$}
  \end{equation*}
$$
\tilde{u} = \tilde{g} 
\quad \text{on $\bar{\mathcal{M}}_{T'}\setminus Q$}, 
$$
where $\tilde{u}(t,x) = u(t,x)$ on $\mathcal{M}_{T'}$,  
$\tilde{u}(T',x) = u(T,x)$,  
$\tilde{g} = g$ on 
$\mathcal{M}_{T'}$ and $\tilde{g}(T',x) = g(T,x)$. 
 Observe that
  \begin{equation*}
    \fddeltatau \tilde{u} = \fddeltatau^{T'} \tilde{u} 
    = \fddeltatau^T u \on \mathcal{M}_{T'}.
  \end{equation*}
  
\end{remark}

\begin{lemma}                       \label{lemma-discrete-comparison-principle}
Assume that 
$a^{\alpha}_k$, 
$b^{\alpha}_k$ and $c^{\alpha}$ satisfy 
Assumptions \ref{assumption 1.23.12.08} and 
\ref{ass-new-for-3}.  
Let $f_1^\alpha$ and $f_2^\alpha$ be functions 
on $A \times \mathcal{M}_{T}$ such that 

$$
\sup_{\alpha} m^\alpha f_2^\alpha 
< \infty, \quad f_1^\alpha \leq f_2^\alpha
\quad \text{on $Q$ for every $\alpha\in A$}.
$$
  Let $u_1$ and $u_2$ be functions 
  on $\bar{\mathcal{M}}_T$ such that 
for some constants 
  $\mu\geq 0$ and  $C\geq 0$ the functions 
  $u_1(t,x)e^{-\mu|x|}$ and  $u_2(t,x) e^{-\mu|x|}$
  are bounded 
on $\bar{\mathcal M}_T$ and
$$
\sup_{\alpha} m^\alpha 
\left(\fddeltatau^T u_1 
+ \Lah u_1 + f_1^\alpha + C\right) 
$$
\begin{equation}                             \label{eqn-comp-princ-ass1} 
\geq \sup_{\alpha} m^\alpha 
\left(\fddeltatau^T u_2 
+ \Lah u_2 + f_2^\alpha \right) \on Q, 
\end{equation}
\begin{equation}                             \label{4.22.12.08}
u_1 \leq u_2\quad\text{on $\bmt \setminus Q$}. 
\end{equation}
  Assume also that $h\leq1$.
  Then there exists a constant
  $\tau^*$ depending only on $K, d_1,
  \mu$ such for $\tau \in (0, \tau^*)$
  \begin{equation}                     \label{eqn-comparison-princ-conclusion}
    u_1 \leq u_2 + T'C \on \bar{\mathcal{M}}_T.
  \end{equation}
  If $u_1$, $u_2$ are bounded on $Q$ then
  \eqref{eqn-comparison-princ-conclusion}
  holds for all positive $\tau$ and $h$.
\end{lemma}

\begin{proof} 
  By using Remark \ref{remark-T'} we may 
  assume that $T=T'$ and 
  $\fddeltatau^T = \fddeltatau$.
Let $w = u_1 - u_2 -C(T'-t)$. Then from
  \eqref{eqn-comp-princ-ass1}
$$
\sup_\alpha m^\alpha \left(\fddeltatau w 
    + \Lah w  \right)\geq
  0, \on Q.
$$
  Notice that, as in \eqref{eq:G2} with $\gamma = 1$ 
  (hence $\xi = 1$ and $\nu = 0$)
  and $f^\alpha = 0$, we have
  $$
G[w] = w + \e \sup_{\alpha}
    m^\alpha \left( \fddeltatau^T w 
      + \Lah w\right)
= \sup_{\alpha \in A}[p^\alpha_\tau \T_\tau w 
  + p^\alpha_k \T_{h,l_k}w + p^\alpha w],
$$
  with
$$
p^\alpha_\tau = \e \tau^{-1} m^\alpha \geq 0, 
  \quad p^\alpha_k = \e m^\alpha (2h^{-2}a^\alpha_h
  + h^{-1}b^\alpha_k) \geq 0
$$
  and 
$$
p^\alpha = 1 - p^\alpha_\tau - \sum_k p^\alpha_k 
  - \e m^\alpha c^\alpha, 
$$
  where one can see that also $p^\alpha \geq 0$
  if $\e$ is sufficiently small. 
  Thus $G$ is a monotone operator 
in the sense that 
  for any $\psi \geq w$ on $\bar{\mathcal M}_T$ we have 
  $G[\psi] \geq G[w]$ on $\mathcal M_T$. 
  So for any sufficiently small fixed $\e > 0$ 
and $\psi\geq w$ on $\bar{\mathcal M}_T$ 
  \begin{equation}                            \label{eqn-w-greater-Gpsi}
    \psi + \e \sup_\alpha m^\alpha 
    \left(\fddeltatau \psi
      + \Lah \psi  \right) \geq w, \on Q.
  \end{equation}
  Let $\gamma \in (0,1)$. Use $\xi$ from
  the proof of Theorem
  \ref{theorem 1.20.3.9}. 
  Then
  \begin{equation*}
    \fddeltatau^T \xi 
= \xi \tfrac{1}{\tau}(\gamma - 1).
  \end{equation*}
  Let $\eta(x) = \cosh(\mu|x|)$  and
  $\zeta = \eta \xi$. Introduce
$$
N_0 = \sup_{\bar{\mathcal M}_T}\tfrac{w_+}{\zeta}.
$$
  Due to the assumption that 
  $u_1(t,x)e^{-\mu|x|}$ and $u_2(t,x)e^{-\mu|x|}$  
  are bounded on 
$\bar{\mathcal M}_T$, we have $N_0 < \infty$.
  Our aim now is to show that, in fact 
  $N_0 = 0$. 

  By Remark
  \ref{remark-taylors-thm},  taking into account that for 
  every $\mu>0$ and integer $n\geq1$ there is a constant $N$ such that 
  for all $x\in\mathbb R^d$  
  \begin{equation*}
    |D_x^n\cosh(\mu|x|)| \leq N \cosh(\mu|x|),
  \end{equation*}
  we get
  $$m^\alpha \Lah \eta(x) \leq 
  m^\alpha L^\alpha \eta(x) 
  + N_1(h^2 + h)\cosh(\mu|x|+\mu K)$$
  \begin{equation}                                          \label{2.1.6.8}
    \leq N_2 \cosh(\mu|x| +\mu K)\leq  N_3 \cosh(\mu|x|),
  \end{equation}
  where $N_1$ and $N_2$ are constants depending on
  $\mu, d_1, K$, and
  $$N_3 := N_2 \sup_{x\in\mathbb R^d}
  \tfrac{\cosh(\mu |x| + \mu K)}{\cosh(\mu|x|)} <
  \infty.$$
  Thus 
  \begin{equation*}
    m^\alpha \left( \fddeltatau
      \zeta + \Lah \zeta \right) 
    \leq \zeta[\tau^{-1}(\gamma-1)+N_3], 
  \end{equation*}
  Let 
  $$\psi := N_0\zeta \geq \zeta
  \tfrac{w_+}{\zeta} \geq w.$$
  Then by \eqref{eqn-w-greater-Gpsi}
  \begin{equation}                                           \label{1.1.6.8}
    w \leq \psi 
    + \e \sup_\alpha 
m^\alpha( \fddeltatau^T \psi + \Lah \psi) 
    \leq \zeta(N_0 + \e\kappa) 
  \end{equation}
  holds on $Q$,  where $\kappa=\kappa(\gamma) =
  \tau^{-1}(\gamma - 1) + N_3$.
  Notice that $\kappa(0) < 0$ for
  $\tau < \tau^* := N_3^{-1}$,   
  and $\kappa(1) >0$. So there is a
  $\gamma\in(0,1)$, which we 
  choose now,  such that
  $\kappa < 0$ and $N_0+\e \kappa > 0$.
  Thus by \eqref{1.1.6.8} and 
\eqref{4.22.12.08}
$$
  w \leq \zeta(N_0 + \e\kappa) \on \bar{\mathcal{M}}_T. 
$$
  Hence
  \begin{equation}                    \label{eqn-proof-of-comp-princ-N0-is-0}
    N_0 = \sup_{\bar{\mathcal{M}}_T} 
\tfrac{w_+}{\zeta} 
    \leq N_0+\e\kappa, 
  \end{equation}
  which implies $N_0=0$, since  $\e\kappa < 0$. 
  This completes the proof of the first assertion
  of the lemma. 

  Assume now that $u_1$ and $u_2$ are bounded on $Q$.
  Then we can take $\mu=0$, i.e., $\eta = 1$. We do 
not need estimate  \eqref{2.1.6.8}, 
hence there is no restriction 
on $h$. We can take $N_3 = 0$ and hence we do not 
  need any restriction on $\tau$.
\end{proof}

\begin{corollary}                                              \label{corollary 1.25.2.9}
Let Assumptions \ref{assumption 1.23.12.08} 
through
\ref{assumption 2.23.12.08} hold. Let $Q$ be a subset of 
$\bar{\mathcal M}_T$. Assume that $g$ is a bounded
function  on $\bar{\mathcal M}_T$ and let $u$ and $w$ denote 
the unique bounded solutions of 
\eqref{eq:fd:a}-\eqref{eq:fd:b} and 
\eqref{1.17.02.09}-\eqref{2.17.02.09}, respectively. 
 Let $\psi$ be a function on 
$\bar{\mathcal M}_T$ such that for some constant $\mu\geq0$ 
the function $e^{-\mu x}\psi(t,x)$ 
is bounded on $\bar{\mathcal M}_T$. Then the following 
statements hold:  
\begin{enumerate}
\item[(i)] Let 
$$
\fddeltatau^T \psi
+ \Lah \psi + f^\alpha \leq 0\quad 
\text{on $Q$ for each $\alpha\in A$}. 
$$
Then $\psi\geq g$ on $\bar{\mathcal M}_T\setminus Q$ 
implies $\psi\geq u$ on $\bar{\mathcal M}_T$, and 
$\psi\geq g$ on $\bar{\mathcal M}_T$ implies  
$\psi\geq w$ on $\bar{\mathcal M}_T$. 
\item[(ii)] Let 
$$
\fddeltatau^T \psi
+ \Lah \psi + f^\alpha \geq 0\quad 
\text{on $Q$ for some $\alpha\in A$}, 
$$
and  $\psi\geq g$ on $\bar{\mathcal M}_T\setminus Q$.  
Then $\psi\leq u$ and $\psi\leq w$ 
on $\bar{\mathcal M}_T$. 
\end{enumerate}
\end{corollary}
\begin{proof}
The statements concerning $u$ follow immediately 
from the previous lemma. Hence the statements 
concerning $w$ follow by 
Remark \ref{remark 3.17.02.09}. 
\end{proof}

Let us consider now problem \eqref{2.30.12.08}-\eqref{3.30.12.08}
and 
\begin{equation}                                              \label{1.18.02.09}
 \max\left[\sup_{\alpha \in A} 
m^\alpha (\fddeltatau^T w 
+ \Lah w + f^\alpha ),g-w\right]  
=  0 \quad \textrm{on}\quad H_T,
\end{equation}
\begin{equation}                                               \label{2.18.02.09}
  w(T,x) =  g(T,x) \quad\text{for $x\in\mathbb R^d$.}                    
\end{equation}

\begin{corollary}                             
                              \label{corollary-soln-to-disc-bellman-pde-is-bdd}

Let Assumptions \ref{assumption 1.23.12.08} 
through
\ref{assumption 2.23.12.08} hold. 
Let $c_1\geq0$ be a constant such that 
\begin{equation}                               \label{1.24.12.08}
\tau^{-1}(e^{c_1\tau}-1)\leq \lambda .
\end{equation}
Then problem \eqref{2.30.12.08}-\eqref{3.30.12.08} 
has a unique bounded solution $u$, and 
\begin{equation}                                \label{2.24.12.08}
  |u(t,x)| \leq N^{\ast}
+e^{-c_1(T-t)}\sup_{x\in\mathbb R^d}|g(T,x)|, 
\end{equation}
holds on $\bar H_T$, where 
\begin{equation}                                  \label{3.18.02.09}
N^{*}=\left\{
      \begin{array}{lcl}
K\rho^{-1}(\lambda^{-1}(1-e^{-\lambda T^{\prime}})+1)
& \text{when} & \lambda> 0,\\
K\rho^{-1}(T^{\prime}+1) &\text{when} & \lambda= 0.
\end{array}\right.
\end{equation}
In addition to the above conditions assume that 
$|g|\leq K$ on $H_{T}$. Then problem 
\eqref{1.18.02.09}-\eqref{2.18.02.09} 
has a unique bounded solution $w$ and 
\eqref{2.24.12.08} holds for $w$ in place of $u$, with 
$2K+1$ in place of $K$ in \eqref{3.18.02.09}.
\end{corollary}
\begin{proof} 
It suffices to prove the corollary for problem 
\eqref{2.30.12.08}-\eqref{3.30.12.08}. Hence we 
get the statement of the corollary 
also for  \eqref{1.18.02.09}-\eqref{2.18.02.09},  
by rewriting it into the form of 
\eqref{2.30.12.08}-\eqref{3.30.12.08}, 
as it is explained in Remark \ref{remark 3.17.02.09}. 
By Theorem \ref{theorem 1.20.3.9} 
problem \eqref{2.30.12.08}-\eqref{3.30.12.08} has a unique 
solution $u$, which is bounded on each grid 
defined by \eqref{4.30.12.08}. Hence it suffices  
to prove estimate 
\eqref{2.24.12.08} on the grid 
$\bar{\mathcal M}_T$. As before, by virtue of Remark
\ref{remark-T'}
  we may assume that $T=T'$ and so $\fddeltatau^T = \fddeltatau$.
Let $\lambda>0$. Then set $N :=\sup_x|g(T,x)|$ and  
  $$
  \xi(t) 
= K\rho^{-1}\{\lambda^{-1}(1-e^{-\lambda(T-t)})+1\} 
+e^{-c_1(T-t)}N.
  $$
  Then on $\mathcal M_T$ 
  $$
  m^\alpha(\fddeltatau \xi + \Lah \xi 
    + f^\alpha)
=m^\alpha\{\fddeltatau \xi 
-\lambda\xi-(c^{\alpha}-\lambda)\xi+f^{\alpha}\}
  $$
  $$  
=-K\rho^{-1}m^{\alpha}
\left[e^{\lambda t-\lambda T}\big(
\frac{e^{\lambda\tau}-1}{\lambda\tau}-1\big)+1\right]  
-m^\alpha (c^{\alpha}-\lambda)\xi
$$
$$
+m^{\alpha}N \tau^{-1}(e^{c_1\tau}-1)e^{c_1(T-t)}-
m^{\alpha}\lambda N e^{-c_1(T-t)}
+m^\alpha f^\alpha.     
  $$
Thus, due to  
$$
\frac{e^{\lambda\tau}-1}{\lambda\tau}>1, 
\quad \xi\geq K\rho^{-1}, 
\quad m^{\alpha}f^{\alpha}\leq K
$$
and conditions \eqref{eqn-ass-m-c} and \eqref{1.24.12.08} we have 
$$
m^{\alpha}(\fddeltatau \xi + \Lah \xi 
    + f^\alpha)\leq 
-K\rho^{-1}m^{\alpha}(1+c^{\alpha}-\lambda)+K
\leq 0\quad\text{on $\mathcal M_T$}. 
$$
Clearly 
$$
\xi(T)\geq \sup_x|g(T,x)| \geq g(T,x).
$$
Hence applying 
Lemma \ref{lemma-discrete-comparison-principle}
with
$u$ and $\xi$ in place of $u_1$ and $u_2$, 
respectively, 
we get $u\leq \xi$ on $\bar{\mathcal{M}}_T$. 
Similarly, by using $-\xi$ in place of $\xi$, 
we get $u\geq-\xi$ 
on $\bar{\mathcal{M}}_T$. If $\lambda=0$ 
then $c_1=0$, and taking $\xi=K\rho^{-1}(T+1)+N$ we 
get \eqref{2.24.12.08} in the same way as above.   
\end{proof}

Finally we can show that Lemma 3.8 
of \cite{krylov:rate:lipschitz} remains valid 
in our setting. 
\begin{lemma} \label{lemma-3.8} 
Assume that Assumptions \ref{assumption 1.23.12.08} 
through \ref{assumption 2.23.12.08} hold.
Let 
$u$ be the solution of
  \eqref{eq:fd:a}-\eqref{eq:fd:b} 
for a bounded function $g$ 
on $\mathbb R^d$. For every integer 
  $n\geq1$ let 
  $f_n^\alpha$ and $g_n$ be functions 
on $A\times H_T$ and on $\mathbb R^d$,  
respectively such that 
$$
\sup_{\alpha\in A}\sup_{\bar H_T}|m^{\alpha}f_n^{\alpha}|
+\sup_{\mathbb R^d}|g_n|\leq K 
\quad\text{for all $n\geq1$}, 
$$  
$$
\lim_{n \to \infty}(\sup_{\alpha} 
m^{\alpha}
|f^\alpha - f^\alpha_n| + |g - g_n|) = 0 
\quad\text{for every $t\in[0,T]$, $x\in\mathbb R^d$}.
$$
  Then $u_n \to u$ 
on $\bar{\mathcal{M}}_T$ as $n \to \infty$, 
where $u_n$ is the 
bounded solution of \eqref{eq:fd:a}-\eqref{eq:fd:b} 
with $f_n^{\alpha}$ and $g_n$ in place of 
$f^{\alpha}$ and $g$, respectively. 
\end{lemma}
\begin{proof}
Having Theorem \ref{theorem 1.20.3.9} 
and Corollary \ref{corollary-soln-to-disc-bellman-pde-is-bdd} 
at our disposal 
we can get 
this lemma in the same way as
Lemma 3.8 in \cite{krylov:rate:lipschitz} is proved: 
Since by Corollary \ref{corollary-soln-to-disc-bellman-pde-is-bdd} 
$u_n$ is bounded uniformly in $n$, any subsequence of 
$\{u_n\}$ contains a subsequence converging to a 
solution of \eqref{eq:fd:a}-\eqref{eq:fd:b}, which is unique 
and equals $u$. Therefore the whole sequence $u_n$ 
converges to $u$. 
\end{proof}

\section{Gradient Estimates 
for Finite Difference Schemes}              \label{section-gradient-estimate}

Thorough this section we assume that Assumption 
\ref{assumption 1.23.12.08} holds.  
Recall that $T'$ denotes the smallest integer
multiple of $\tau$ which is greater than
or equal to $T$. For a fixed number  
$\e\in (0,Kh]$ and a unit vector $l \in R^d$, 
set $h_r =h$ for $r = \pm 1, \ldots, \pm d_1$ and
$h_r = \e$ for $r = \pm (d_1 + 1)$, and 
$\ell_{\pm(d_1+1)}=\pm l$. 
Define 
\begin{equation*}
  \bar{\mathcal{M}}_T(\e) := \{(t,x+i\e l):(t,x) 
\in \bar{\mathcal{M}}_T, i = 0, \pm 1, \ldots  \},
\end{equation*}
\begin{equation*}
  {\mathcal{M}}_T(\e) :=  \bar{\mathcal{M}}_T(\e) 
\cap([0,T)\times\mathbb R^d). 
\end{equation*}
Let $Q \subset
\bar{\mathcal{M}}_T(\e)$ be a nonempty finite set. Define  
$Q^{\prime}=Q\cap([0,T)\times \mathbb R^d$),
\begin{gather*}
  Q^0_\e = \{(t,x): (t+\taut, x) \in 
Q^{\prime}, (t,x+h_r\ell_r) \in Q, 
 \forall r = \pm 
  1, \ldots , \pm (d_1 + 1)\}
\end{gather*}
and $\partial_\e Q := Q \setminus
Q^0_\e$. 

\begin{assumptions}                        \label{ass-for-thm-5.2} 
  For  $r =\pm 1, \ldots, \pm (d_1 + 1)$ 
and  
$\alpha\in A$ on $Q^0_\e$  we have 
  \begin{equation}                         \label{eqn-first-ass-for-thm-5.2}
    |\fddelta_{h_r,\ell_r} b^\alpha_k| 
\leq K, 
\quad 
m^\alpha|\fddelta_{h_r,\ell_r} f^\alpha| \leq K, 
\quad
m^\alpha|\fddelta_{h_r,\ell_r} c^\alpha|  \leq K,
  \end{equation}
  \begin{equation}                         \label{eqn-second-ass-for-thm-5.2}
    |\fddelta_{h_r,\ell_r} a^\alpha_k|
\leq K\sqrt{a^\alpha_k} + K h.
  \end{equation}
\end{assumptions}

The following estimate plays a crucial role 
in the proof of Theorem \ref{theorem-the-main-result}. 
It generalizes Theorem 5.2 from 
\cite{krylov:rate:lipschitz:published}. 

\begin{theorem}                                \label{thm-est-disc-grad}
  Let 
Assumptions
 \ref{assumption 1.23.12.08}, 
\ref{assumption 2.23.12.08} 
and 
\ref{ass-for-thm-5.2}
hold. 
Let $u$ be a function on
$\bar{\mathcal{M}}_T(\e)$ such that it 
satisfies 
(\ref{eq:fd:a}) with $Q^{\prime}$ 
in place of $Q$.
Then there is a constant $N^{\ast}\geq0$, depending 
only on $d_1$ and  $K$  
such that for any constant $c_o\geq 0$ satisfying
  \begin{equation}                            \label{assumption-on-c_0}
    \lambda + \tfrac{1}{\tau}(1-e^{-c_o \tau}) 
\geq N^{\ast}+1,
  \end{equation}
we have 
\begin{equation}                             \label{eqn-est-dis-grad}
    |\delta_{\e,\pm l}u| 
\leq \sqrt{\tfrac{2N^{\ast}}{\rho}}
e^{c_o (T+\tau)}
\big(1+\max_{Q}|u| 
+ \max_{r=\pm1,\dots,\pm(d_1+1)}\max_{\partial_\e Q}
|\delta_{h_r,\ell_r}u|)
\quad\text{on $Q$}. 
\end{equation}
In addition to the above conditions assume that 
$g$ is a function on $\bar H_T$ such that 
$|\fddelta_{h_r,l_r}g|\leq K$ on $Q^{0}_{\varepsilon}$ for every 
$r=\pm1,\dots,\pm(d_1+1)$.  
Let $w$ be a
function on 
$\bar{\mathcal{M}}_T(\e)$ that 
satisfies 
\eqref{1.17.02.09} with $Q^{\prime}$ 
in place of $Q$. Then the above statement holds also 
for $w$ in place of $u$. 
\end{theorem}

\begin{proof}
We follow the proof of Theorem 5.2 from 
\cite{krylov:rate:lipschitz:published} with some changes. 
Let
\begin{equation*}
v_r = \fddelta_{h_r,l_r}v, \quad v = \xi u, \quad \xi(t) 
=  \left\{
      \begin{array}{ll}
        e^{c_o t} & t < T, \\
        e^{c_o T'} & t = T
      \end{array} \right\},   
\end{equation*}
where $T'$ denotes the smallest multiple of $\tau$ that is not less than $T$. 
Let $(t_0,x_0) \in Q$ be the point where
$$V = \sum_r (v_r^-)^2$$
is maximized. By definition, 
for any $(t,x) \in Q_\varepsilon^o$ we know that 
$$(t,x+h_r\ell_r) \in Q.$$
Clearly, either
$$v_r(t,x) \leq 0 \quad \textrm{or} 
\quad -v_r (t,x) = v_{-r}(t,x+h_r\ell_r) 
\leq 0.$$
Consequently, 
$$
|v_r(t,x)| \leq V^{1/2}(t_0,x_0).
$$
Hence
\begin{gather}
M_1 := \sup_{Q,r}|v_r| 
\leq \sup_{\partial_\varepsilon Q,r}|v_r| 
+
V^{1/2}(t_0,x_0),                               \label{5.5} \\ 
|\delta_{\varepsilon,
\pm l}u| 
\leq e^{c_0 T'} 
\sup_{\partial_\varepsilon Q,r}|\delta_{h_r,\ell_r}u| 
+ V^{1/2}(t_0,x_0)
\quad \text{on $Q$}.
\end{gather}
So we need only estimate $V$ on $Q$. 
If $(t_0,x_0)$ belongs to $\partial_{\varepsilon} Q$, 
then the conclusion of the
theorem is clearly true. Thus, we may assume that $(t_0,x_0) \in
Q^0_\varepsilon$. For any ${\varepsilon}_0 > 0$ 
there exists $\alphao \in
A$ such that at $(t_0,x_0)$, 
\begin{equation*}
m^\alphao \left(\fddeltatau^T u 
+ a^{\alphao}_k \fdDeltak u + b^{\alphao}_k
\fddeltak u - c^\alphao u + f^{\alphao} \right) 
+ {\varepsilon}_0 \geq 0, 
\end{equation*}
and so for some ${\varepsilon} 
\in [0,{\varepsilon}_0]$
\begin{equation}\label{5.8}
m^\alphao \left(\fddeltatau^T u 
+ a^{\alphao}_k \fdDeltak u + b^{\alphao}_k
\fddeltak  u - c^\alphao u + f^{\alphao} \right) 
+ {\varepsilon} = 0.
\end{equation}
Furthermore (thanks to the fact that $(t_0,x_0) 
\in Q^0_\varepsilon$)
\begin{equation}\label{5.9}
\Thrlr \left[m^\alphao \left(\fddeltatau^T u + a^{\alphao}_k 
\fdDeltak u + b^{\alphao}_k \fddeltak u  
- c^\alphao u + f^{\alphao} \right) 
\right] \leq 0, 
\end{equation}
where $\T_{h,l}\varphi(t,x):=\varphi(t,x+hl)$ 
for any number $h$, 
vector $l\in\mathbb R^d$ and function $\varphi$ defined 
at $(t,x)$ and $(t,x+hl)$.  
Here and below $(t_0,x_0)$ is fixed and 
for simplicity of notation it 
is omitted in the arguments of the functions.
We subtract
(\ref{5.8}) from (\ref{5.9}) and divide by
$h_r$ to obtain that for each $r$
\begin{equation*}
 m^{\alphao} \fddelta_{h_r,\ell_r}
\left(\fddeltatau^T u + a^{\alphao}_k 
\fdDeltak u + b^{\alphao}_k \fddeltak u 
+ f^{\alphao} - c^\alphao u\right)  -
\tfrac{{\varepsilon}}{h_r} \leq 0.
\end{equation*} 
By the discrete Leibnitz rule 
\begin{equation}\label{5.13}
\begin{split}
m^\alphao \left(\delta_\tau (\xi^{-1} v_r) + \xi^{-1}
\left[a_k^{\alphao} 
\fdDeltak v_r +I_{1r} + I_{2r} + I_{3r}
  \right] + \fddelta_{h_r,l_r} f^\alphao \right)        \\ 
- \xi^{-1}\fddelta_{h_r,\ell_r}(m^\alphao c^\alphao v)  -
\tfrac{{\varepsilon}}{h_r} \leq 0,
\end{split}
\end{equation}
where
\begin{eqnarray*}
I_{1r} & = & ( \fddelta_{h_r,\ell_r} a^{\alphao}_k ) 
\fdDeltak v,                                              \\
I_{2r} & = & h_r ( \fddelta_{h_r,\ell_r} 
a^{\alphao}_k) \fdDeltak v_r,                             \\
I_{3r} & = & 
b^{\alphao}_k
\fddeltak v_r + (
\fddelta_{h_r,\ell_r} b^{\alphao}_k )  
\Thrlr\fddeltak v.  
\end{eqnarray*}
Notice that
\begin{eqnarray*}
0 & \geq & \fdDeltak \sum_r (v_r^-)^2 = 2v_r^- 
\fdDeltak v_r^- + \sum_r
\left[ 
(\fddeltak v_r^-)^2 
+ 
(\fddelta_{h_k,\ell_{-k}}
v_r^-)^2\right]                                               \\
& \geq & - 2v_r^- \fdDeltak v_r 
+ \sum_r (\fddeltak v_r^-)^2 
+ \sum_r (\fddeltamk v_r^-)^2,
\end{eqnarray*}
which gives 
\begin{equation}                               \label{22.12.08}
0 \leq v_r^- \fdDeltak v_r
\end{equation}
and
\begin{equation*}
I := \sum_r a^{\alphao}_k(\fddeltak v_r^-)^2 
\leq v_r^-
a^{\alphao}_k\fdDeltak v_r. 
\end{equation*}
Multiplying \eqref{5.13}    
by $\xi v_r^-$ 
 and summing 
up in $r$ we get
\begin{equation}\label{5.13.1}
\begin{split}
m^\alphao 
\Big(
\xi v_r^- \delta_\tau (\xi^{-1} v_r) 
+ \tfrac{1}{2}
a^{\alphao}_k v_r^- \fdDeltak v_r 
+ \tfrac{1}{2}I  + v_r^-
[I_{1r} + I_{2r} \\ + I_{3r} 
+ \xi \fddelta_{h_r,\ell_r} f^\alphao 
  ] 
\Big) 
- v_r^- \fddelta_{h_r,\ell_r}(m^\alphao c^\alphao v)
 - \xi 
v_r^- 
\tfrac{{\varepsilon}}{h_r} \leq 0.
\end{split}
\end{equation}
Since $-v_r^-v_r = \sum_r(v_r^-)^2$, 
$m^\alphao \fddelta_{h_r,\ell_r} f^\alphao
\geq -K$ and  
$m^{\alpha_0}|\fddelta_{h_r,\ell_r} c^{\alphao}|\leq K$, 
we have 
\begin{multline*}
m^\alphao v_r^- 
\xi \fddelta_{h_r,\ell_r} f^\alphao - v_r^- 
\fddelta_{h_r,\ell_r}(m^\alphao c^\alphao v)              \\ 
= m^\alphao v_r^- 
\xi \fddelta_{h_r,\ell_r} f^\alphao - m^\alphao\vrm
(\fddelta_{h_r,l_r} c^{\alphao}) 
\Thrlr v  
- m^\alphao\vrm 
c^{\alphao} v_r                           \\
\geq - e^{c_0T'}
 K \sum_{r}\vrm  -m^{\alphao} 
\vrm |\fddelta_{h_r,\ell_r}
c^{\alphao}| |\Thrlr v| 
+ m^\alphao
c^{\alphao}
\sum_r(\vrm)^2                                    \\
 \geq - e^{c_o T'}2K(d_1+1)M_1  - 
2(d_1+1)K M_1 M_0 
+  m^{\alpha_0}c^{\alpha_0} V,
\end{multline*}
where 
$$
M_0:=\max_{Q}|v|. 
$$
Since  $V$ attains its maximum at 
$(t_0,x_0)\in Q_{\varepsilon}^0$ we have 
\begin{eqnarray*}
0 & \geq & \sum_r \fddeltak (\vrm)^2 = 2\vrm
\fddeltak \vrm + \sum_r h_k (\fddeltak \vrm)^2 \\
& \geq & 2 \vrm \fddeltak \vrm \geq -2 \vrm \fddeltak v_r. 
\end{eqnarray*}
Next recall that $b_k^\alpha \geq 0$ and
$|\fddelta_{h_r,l_r}b^\alpha_k| \leq K$. 
Therefore
$$
-\vrm b^{\alphao}_k
\fddeltak v_r\leq 0,
$$
and 
$$
v_r^- I_{3r}  \geq  -\vrm |\fddelta_{h_r,\ell_r}
b^{\alphao}_k||\Thrlr \fddeltak v| 
\geq -4Kd_1(d_1+1)M_1^2.
$$
By the discrete Leibnitz rule 
\begin{eqnarray*}
\xi \vrm \fddeltatau^T (\xi^{-1} v_r ) 
& = & 
\xi \vrm \left[\xi^{-1}(t_0
  +\tau )\fddeltatau^T v_r + v_r \fddeltatau^T \xi^{-1}  \right] \\
& = & 
e^{-c_o\tau} \vrm \fddeltatau^T v_r - V \xi \fddeltatau^T \xi^{-1} 
\geq
- V\xi \fddeltatau^T \xi^{-1}                                     \\
& = &\nu V, 
\end{eqnarray*}
where 
$$
\nu=\nu(c_o)=\tfrac{1}{\tau}(1-e^{-c_o\tau}).
$$
Using the above estimates we get 
\begin{equation*}
\begin{split}
m^\alphao (\nu +c^\alphao) V 
+ \tfrac{1}{2} m^{\alphao} a^{\alphao}_k \vrm
\fdDeltak v_ r +
\tfrac{1}{2} m^{\alphao} I                                         \\ 
+ \vrm m^\alphao [I_{1r} + I_{2r}] - \xi \vrm
\tfrac{{\varepsilon}}{h_r} 
\leq 
2(d_1+1)KM_1(e^{c_oT'} + M_0 + 2d_1m^\alphao M_1).  
\end{split}
\end{equation*}
Hence
\begin{equation*}
\begin{split}
m^\alphao (\nu+c^\alphao)V 
\leq 2(d_1+1)KM_1(e^{c_oT'} + M_0 + 2d_1m^\alphao M_1)   \\
 +m^\alphao
\vrm|I_{1r}| + m^\alphao \vrm|I_{2r}| - \tfrac{1}{2}m^\alphao
a^{\alphao}_k \vrm \fdDeltak v_r - 
\tfrac{1}{2}m^\alphao I + \xi \vrm
\tfrac{{\varepsilon}}{h_r}.
\end{split}
\end{equation*}
Define
\begin{align*}
J_1 & := \vrm |(\fddelta_{h_r,\ell_r} a^{\alpha_o}_k) \fdDeltak v| -
\tfrac{1}{4}\sum_r a^{\alphao}_k (\fddeltak \vrm)^2,\\
J_2 & := J_3 - \tfrac{1}{2} a^{\alphao}_k \vrm \fdDeltak v_r -
\tfrac{1}{4}\sum_r a^{\alphao}_k (\fddeltak \vrm)^2,
\end{align*}
$$
J_3 := h_r \vrm |(\fddelta_{h_r,\ell_r} a^{\alphao}_k) \fdDeltak
v_r|.
$$
Then we can rewrite the above inequality as
\begin{equation}                                      \label{2.22.12.08}
m^\alphao (\nu+c^\alphao) V 
\leq 
2(d_1+1)KM_1(e^{c_oT'} + M_0 
+ 2d_1m^\alphao M_1) 
\end{equation}
\begin{equation*}
+ m^\alphao(J_1 + J_2) 
+ \xi \vrm
\tfrac{{\varepsilon}}{h_r} .
\end{equation*}
So we need to estimate $J_1, J_2$.
We turn our attention to $J_1$. 
Using 
condition \eqref{eqn-second-ass-for-thm-5.2},  
noticing  
that $h|\fdDeltak v|\leq 2M_1$ and 
\begin{equation*}
|\fdDeltak v| \leq \sum_r |\fddeltak \vrm| 
+ \sum_r |\fddelta_{h_k,\ell_{-k}} \vrm| , 
\end{equation*}
we have  
\begin{eqnarray*}
\vrm |(\fddelta_{h_r,l_r} a_k^\alphao) \fdDeltak v| 
& 
\leq & N_1
M_1 |(\sqrt{a_k^\alphao} + h) \fdDeltak
v|                                                              \\ 
& 
\leq 
& 
N_1
M_1\sqrt{a_k^\alphao} |\fdDeltak v| + 
N_2M_1^2                                         \\
& \leq & 
2N_1M_1 \sqrt{a_k^\alphao} 
\sum_r |\fddeltak\vrm|
+ 
N_2M_1^2                                         \\ 
& \leq & 
N_3
M_1^2
+ \tfrac{1}{4} 
\sum_r a^\alphao_k (\fddeltak \vrm)^2, 
\end{eqnarray*}
where $N_1$, $N_2$ and $N_3$ are constants 
depending only on $d_1$ and $K$. 
So 
\begin{equation}                                 \label{3.22.12.08}
J_1 \leq N_3M_1^2.
\end{equation}
Next we estimate $J_3$. 
Since $h_r \leq Kh$ for all $r$, 
by condition \eqref{eqn-second-ass-for-thm-5.2}
\begin{equation*}
  J_3 \leq K^2
hv_r^- \sqrt{a_k^\alphao} | \fdDeltak v_r| 
+ K^2h^2 \sum_k\vrm |\fdDeltak v_r|.
\end{equation*}
Hence using  $h^2 |\fdDeltak v_r| \leq 4M_1$
 and $|a| = 2a^{-}+a$,  we get
\begin{equation*}
  J_3 \leq 2K^2h 
\vrm \sqrt{a_k^{\alphao}} (\fdDeltak v_r)^{-} 
+ K^2h \vrm \sqrt{\alpha_k^\alphao} \fdDeltak v_r + 
8K^2d_1M_1^2.
\end{equation*}
Notice that the summations 
in $r$ above can be restricted to $\{r:v_r<0\}$. 
For these $r$  we have 
\begin{equation*}
  h (\fdDeltak v_r)^{-}  
\leq h|\fdDeltak \vrm| 
\leq |\fddeltak \vrm| + |\fddelta_{h_k,\ell_{-k}} \vrm|. 
\end{equation*}
Hence 
\begin{align*}       
J_3 & \leq 4K^2
\vrm \sqrt{a_k^{\alphao}} 
|\delta_{h_k,\ell_k} v_r^{-}|+K^2\vrm h
\sqrt{a_k^\alphao}
\fdDeltak v_r  + 8K^2d_1M_1^2                                            \\ 
& \leq NM_1^2+ \tfrac{1}{4}\sum_r a_k^\alphao (\fddeltak \vrm)^2 
+K^2\vrm h \sqrt{a_k^\alphao}\fdDeltak v_r,                           \\
  J_2 & \leq NM_1^2  - \tfrac{1}{2}
\big(a_k^\alphao - 2K^2h\sqrt{a_k^\alphao}\big) 
\vrm \fdDeltak v_r.
\end{align*}
By \eqref{22.12.08}
$$
J_2 \leq NM_1^2  - 
\tfrac{1}{2}\sum_{k\in\mathcal K}R_k, 
$$
where 
$$
R_k=(a_k^\alphao - 2K^2h\sqrt{a_k^\alphao})\vrm \fdDeltak v_r, 
\quad
\mathcal{K} := \left\{k : a_k^\alphao - 2K^2h 
\sqrt{a_k^\alphao}<0\right\}.
$$
Notice that for $k\in\mathcal K$ we have 
$a_k^{\alphao}<4K^4h^2$ and hence 
$$
|R_k|\leq 4K^4h^2|v^{-}_r\Delta_{h_k,\ell_k}v_r|\leq NM_1^2
$$
with a constant $N$ depending only on $K$ and $d_1$.
Thus $J_2 \leq NM_1^2$  and  hence 
by \eqref{2.22.12.08}  and 
\eqref{3.22.12.08} we get 
\begin{equation*}
m^\alphao (\nu+c^\alphao)V 
\leq NM_1\left(e^{c_oT'} + M_0 + m^\alphao M_1\right) 
+ \xi \vrm
\tfrac{{\varepsilon}}{h_r} ,
\end{equation*}
where $N$ denotes constants 
depending only on $K$ and $d_1$.
By \eqref{5.5} we have $M_1 \leq \mu + V^{1/2}$, 
where 
$$
\mu := \sup_{\partial_\varepsilon Q,r} |v_r| 
\leq e^{c_0 T'} \bar \mu, 
\quad 
\bar \mu=|\sup_{\partial_\varepsilon Q,r}\fddelta_{h_r,l_r} u|.
$$ 
Set  
$$
\bar{M_0} = |u|_{0,Q} 
\geq e^{-c_0 T'}M_0, \quad \bar{V} = e^{-2c_0 T'}V. 
$$
Then, using Young's inequality,
we obtain 
\begin{multline}                                             \label{eqn-5.15}
m^\alphao (\nu+c^\alphao)\bar V
\leq N(\bar{\mu}+\bar{V}^{1/2})\left(1 
+ \bar{M_0} + \bar{\mu} + m^\alphao \bar{V}^{1/2}\right) 
+  e^{-c_0T^{\prime}}\vrm
\tfrac{\bar{\varepsilon}}{h_r}  \\
\leq N^*\left(1 + \bar{M_0}^2 
+ \bar{\mu}^2 + m^\alphao \bar{V}\right) + 
\frac{\rho}{2}\bar V+
e^{-c_oT^{\prime}}\vrm
\tfrac{{\varepsilon}}{h_r}.
\end{multline}
Assume that for $c_o$
$$
\lambda + \nu(c_o) 
\geq 1 + N^*. 
$$
Then (\ref{eqn-5.15}) yields 
\begin{equation*}
 m^\alphao \left(1+c^\alphao - \lambda \right)
\bar V 
\leq N^*\left(1 + \bar{M_0}^2 + \bar{\mu}^2\right) + 
\frac{\rho}{2}\bar V+
e^{-c_oT^{\prime}}
\vrm\tfrac{{\varepsilon}}{h_r}.
\end{equation*}
Hence using condition \eqref{eqn-ass-m-c}  
and then letting ${\varepsilon} \to 0$  
 we obtain
\begin{equation*}
\bar V 
\leq \frac{2}{\rho}N^*\left(1 + \bar{M_0}^2 + \bar{\mu}^2\right), 
\end{equation*}
that obviously yields estimate \eqref{eqn-est-dis-grad}. 
Finally we use Remark to rewrite 
equation for $w$ into equation 
(\ref{eq:fd:a}) with $Q^{\prime}$ 
in place of $Q$, and notice that for 
$\gamma=(\alpha,r)\in\bar A$
$$
m^{\gamma}|\delta_{h_r,\ell_r}f^{\gamma}|
=m^{\alpha}(1+r)^{-1}
|\delta_{h_r,\ell_r}f^{\alpha}
+\tfrac{r}{m^{\alpha}}\delta_{h_r,\ell_r}g|\leq 2K
$$
for $r=\pm1,\dots,\pm(d_1+1)$. Hence the statement 
on $w$ follows from that on $u$.  
\end{proof}

Let us consider now \eqref{eq:fd:a}-\eqref{eq:fd:b} 
with 
$\mathcal M_T(\varepsilon)$ in place of $Q$.

\begin{corollary}                          \label{col-grad-est}
Assume that 
Assumptions
 \ref{assumption 1.23.12.08} through  
\ref{assumption 2.23.12.08} 
and 
\ref{ass-for-thm-5.2} with $\mathcal M_T(\varepsilon)$ 
in place of $Q^{0}_{\varepsilon}$
hold. Let $g$ be a bounded 
function on $\mathbb R^d$. Let $u$ be the solution to
\eqref{eq:fd:a}-\eqref{eq:fd:b}
with $Q=\mathcal M_T(\e)$.  
Then there is a constant $N^{\ast}\geq0$, 
depending only on $d_1$ and $K$ 
 such that 
for any constant $c_o\geq0$ 
satisfying \eqref{assumption-on-c_0} we have 
\begin{equation}                                   \label{eqn-grad-est-col}
    |\delta_{\e,\pm l}u| 
\leq N_0e^{c_o(T+\tau)}
\big(N_1+\sup_{\mathbb R^d}|g|
+\max_{r}\sup_{\mathbb R^d}|\delta_{h_r,\ell_r}g|\big)
\quad\text{on $\bar{\mathcal M}_T(\e)$,}
  \end{equation}
  where $N_0$ and $N_1$ are constants. 
The constant $N_0$ depends  
only on $K$, $d_1$ and $\rho$ and the constant $N_1$ 
depends on $K$, $d_1$, $\rho$ and $\lambda$, 
provided $\lambda>0$, and if $\lambda=0$ 
then it depends on $K$, $d_1$, $\rho$ and $T$. 
\end{corollary}

\begin{proof}
Let $B_r$ denote the open ball of radius $r$ centered 
at the origin in $\mathbb R^d$. Using Theorem
\ref{thm-est-disc-grad} 
  with $Q_n := \bar{\mathcal M}_T(\e) \cap ([0,T] \times B_n)$ 
in place of $Q$ for any integer 
$n\ge1$ we have
\begin{equation*}
    |\delta_{\e,\pm l}u| \leq Ne^{c_o(T+\tau)}
\big(1+\max_{Q_n}|u| +\max_{r}\max_{\partial_\e Q_n}
|\delta_{h_r,\ell_r} u|\big), \quad \text{on $Q_n$},
\end{equation*}
where $N$ is a constant depending only on $d_1$ and $K$. 
In addition to the assumptions assume that 
for all $\alpha\in A$ the functions $f^{\alpha}$ 
and $g$ vanish  outside of 
a fixed ball of radius $R$ centered at the origin 
in $\mathbb R^d$. Set 
$\partial^T_{\varepsilon}
=\{(T,x)\in\partial_{\varepsilon}Q_n\}$. 
Then by Lemma \ref{corollary-R}
\begin{equation*}
    \lim_{n \to \infty} 
\sup_{k,\partial_\e Q_n \setminus \partial^T_{\varepsilon} Q_n}     
\left(|\fddeltak u| + |\delta_{\e,\pm l}u| \right)
= 0. 
\end{equation*}
Hence on $\bar{\mathcal M}_T(\e)$
  \begin{equation}                                              \label{eqn-pf-col-4.3}
    |\delta_{\e,\pm l}u| 
\leq Ne^{c_oT'}\big(1+\sup_{\bmt(\e)}|u|
+ \max_r\sup_{\mathbb R^d}
|\delta_{h_r,\ell_r}g|\big).
  \end{equation}
Let us now remove the additional assumption on 
$f^\alpha$ and $g$. Let $\eta\in C_0^{\infty}(\mathbb R^d)$ 
be a nonnegative function such that $\eta\leq 1$, 
$|D\eta|\leq1$ on the whole $\mathbb R^d$ 
and $\eta(x)=1$ for $|x|\leq1$. 
For each integer $n\geq1$ define 
$$
f^{\alpha}_n(t,x)=\eta(n^{-1}x)f^{\alpha}(t,x), 
\quad g_n(x)=\eta(n^{-1}x)g(x), 
\quad t\geq0, \,x\in\mathbb R^d .  
$$
Then clearly
  \begin{equation*}
    \lim_{n \to \infty}
(\sup_{\alpha} m^{\alpha}|f^\alpha - f^\alpha_n| + |g - g_n|) = 0
\quad \text{on $\bar H_{T}$}, 
  \end{equation*}
$$
|f^{\alpha}_n|\leq |f^{\alpha}|, 
\quad |\delta_{h_r,\ell_r}f^{\alpha}_n|
\leq
|\ell_r|\sup_{\bar H_{T}}|f^{\alpha}|
+|\delta_{h_r,\ell_r}f^{\alpha}|, 
$$
\begin{equation}                              \label{1.27.12.08}
|g_n|\leq |g|,\,
\quad
|\delta_{h_r,\ell_r}g_n|
\leq
|\ell_r|\sup_{\mathbb R^d}|g|+|\delta_{h_r,\ell_r}g|. 
\end{equation}
 Let $u_n$ be the
  solution to
  \eqref{eq:fd:a}-\eqref{eq:fd:b} 
with $Q=\bar{\mathcal M}_T(\e)$ and with 
$f_n^{\alpha}$ and $g_n$ 
in place of $f^{\alpha}$ and $g$, respectively. 
Then from
  \eqref{eqn-pf-col-4.3} 
and \eqref{1.27.12.08} 
for all $n \in \N$,
  \begin{equation*}
    |\delta_{\e,\pm l}u_n| 
\leq Ne^{c_o(T+\tau)}
\big(1+\sup_{\bar{\mathcal M}_T(\e)}|u_n| + 
K\sup_{\mathbb R^d}|g|
+\max_r\sup_{\mathbb R^d}
|\delta_{h_r,\ell_r}g|\big).
  \end{equation*}
Hence estimating 
$\sup_{\bar{\mathcal M}_T(\e)}|u_n|$ by using
Corollary \ref{corollary-soln-to-disc-bellman-pde-is-bdd} 
and then letting $n\to\infty$ by
using Lemma \ref{lemma-3.8} we get 
estimate \eqref{eqn-grad-est-col}. 
\end{proof}

\begin{assumptions}                        \label{assumption 1.02.01.08}
For all $\alpha\in A$, $t\geq0$ and $x,y\in\mathbb R^d$
$$
|b^\alpha_k(t,x)-b^\alpha_k(t,y)| 
\leq K|x-y|, 
\quad 
m^{\alpha}|c^\alpha(t,x)-c^{\alpha}(t,y)|  \leq K|x-y|, 
$$
$$
m^\alpha|f^\alpha(t,x)-f^\alpha(t,y)| \leq K|x-y|,
$$
\begin{equation}                         \label{assumption 1.30.12.08}
    |\sqrt{a^{\alpha}(t,x)}-\sqrt{a^{\alpha}(t,y)}|
\leq K|x-y|.
 \end{equation}
\end{assumptions}

\begin{theorem}                          \label{theorem 4.17.3.9}
Let Assumptions 
\ref{assumption 1.23.12.08} through 
\ref{assumption 2.23.12.08} 
and Assumption \ref{assumption 1.02.01.08} hold. 
Assume that $g$ is a Borel function on $\bar H_T$  
such that 
$$
\sup_{\bar H_T}|g|\leq K, \quad |g(t,x)-g(ty)|\leq K|x-y|
\quad\text{for all $t\in[0,T]$, $x,y\in\mathbb R^d$}. 
$$ 
Then there is a constant $N^{\ast}\geq0$ such that for any constant 
$c_o\geq0$ satisfying \eqref{assumption-on-c_0} for the solution 
$u$ of \eqref{2.30.12.08}-\eqref{3.30.12.08} 
and the solution $w$ of 
\eqref{1.18.02.09}-\eqref{2.18.02.09} we have 
\begin{equation}                                       \label{1.31.12.08}
    |u(t,x) - u(t,y)| +|w(t,x) - w(t,y)|
\leq Ne^{c_o(T+\tau)}|x-y|, 
\end{equation}
for all $t\in[0,T]$, $x,y\in\mathbb R^d$, 
where $N$ is a constant, that depends only on 
$K$, $d_1$, $\rho$ and $\lambda$, if $\lambda>0$. 
If $\lambda=0$ then $N$ depends on 
$K$, $d_1$, $\rho$ and $T$. 
\end{theorem}

\begin{proof}
To prove \eqref{assumption 1.30.12.08} 
let $(t,x)$ and $(t,y)$ be fixed elements of 
$\bar H_T$. We may assume that $t<T$. Moreover, 
by making a suitable shift in the argument of 
the functions,  
we may assume that $(t,x)\in\mathcal M_T$. 
If $|x-y|\geq K$, then estimate \eqref{1.31.12.08} 
holds by virtue of 
Corollary  
\ref{corollary-soln-to-disc-bellman-pde-is-bdd}.  
Assume that $|x-y|<K$. Set $\ell=(y-x)/|x-y|$, 
$\ell_{\pm(d_1+1)}=\pm\ell$ and 
$\varepsilon=|x-y|/{n}$, where $n$ is the smallest 
positive integer such that $|x-y|/{n}\leq Kh$. 
Then 
$$
 |u(t,x) - u(t,y)| 
\leq\varepsilon
\sum_{j=0}^{n-1}|\delta_{\varepsilon,\ell}u(t,x+j\varepsilon\ell)|
$$
$$
\leq n\varepsilon
\sup_{\bar{\mathcal M}_T(\varepsilon)}|\delta_{\varepsilon,\ell}u|
=|x-y| 
\sup_{\bar{\mathcal M}_T(\varepsilon)}|\delta_{\varepsilon,\ell}u|. 
$$
Hence we can finish the proof 
by using Corollary \ref{col-grad-est} if we show  
that Assumption \ref{ass-for-thm-5.2} with $\mathcal M_T(\varepsilon)$ 
in place of $Q^{0}_{\varepsilon}$ holds. 
It is easy to see that condition \eqref{eqn-first-ass-for-thm-5.2} 
is satisfied with $K^2$ in place of $K$. To verify 
condition \eqref{eqn-second-ass-for-thm-5.2} notice that 
for any $r=\pm1,\dots\pm(d_1+1)$, $\ell_r$ and 
$(t,z)\in\bar{\mathcal M}_T(\varepsilon)$
$$
|\delta_{h_r,\ell_r} a_k^{\alpha}(t,z)|
=h_r^{-1}2(|a_k^{\alpha}(t,z)|^{1/2}-|a_k^{\alpha}(t,z+h_r\ell_r)|^{1/2})
$$
$$
+h_r^{-1}
\big(
|a_k^{\alpha}(t,z+h_r\ell_r)|^{1/2}
-|a_k^{\alpha}(t,z+h_r\ell_r)|^{1/2}
\big)^2
$$
$$
\leq K^2|a_k^{\alpha}(t,z)|^{1/2}+h_rK^4
\leq K^{\prime}|a_k^{\alpha}(t,z)|^{1/2}
+K^{\prime}h
$$
 with $K^{\prime}:=1+K^4$. The proof is complete.
\end{proof}

Now we investigate the dependence of the solution 
to \eqref{eq:fd:a}-\eqref{eq:fd:b} on the data. 
Therefore together with  
$a^{\alpha}_k$, $b^{\alpha}_k$, $c^{\alpha}$, 
$f^{\alpha}$ we consider also functions 
$\hat a^{\alpha}_k$, $\hat b^{\alpha}_k$, 
$\hat c^{\alpha}$, $\hat f^{\alpha}$ 
defined on $H_{\infty}$ for each $\alpha\in A$. 
   
\begin{assumptions}                    \label{assumption 1.01.01.09}
Assumptions 
\ref{assumption 1.23.12.08} through 
\ref{assumption 2.23.12.08} and 
Assumption \ref{ass-for-thm-5.2}  
with $\mathcal M_T(0)=\mathcal M_T$ and 
$r=\pm1,\dots,\pm d_1$ 
in place of $Q^{0}_{\varepsilon}$ 
and $r=\pm1,\dots,\pm(d_1+1)$, respectively, 
hold for $a^{\alpha}_k$, $b^{\alpha}_k$, $c^{\alpha}$ 
and 
$f^{\alpha}$ and also 
for 
$\hat a^{\alpha}_k$, $\hat b^{\alpha}_k$, 
$\hat c^{\alpha}$ 
and $\hat f^{\alpha}$ in place of 
$a^{\alpha}_k$, $b^{\alpha}_k$, $c^{\alpha}$  
and $f^{\alpha}$, respectively, with the same 
function $m^{\alpha}$ and constant $\lambda\geq0$.  
\end{assumptions}

If Assumption \ref{assumption 1.01.01.09} 
holds and $g$ and $\hat g$ are bounded functions 
on $\mathbb R^d$, then by 
Theorem \ref{theorem 1.20.3.9}
we have, in particular, the existence of a unique bounded solution  
of \eqref{2.30.12.08}-\eqref{3.30.12.08} with 
$a^{\alpha}_k$, $b^{\alpha}_k$, $c^{\alpha}$, 
$f^{\alpha}$ and $g$ and also with  
$\hat a^{\alpha}_k$, $\hat b^{\alpha}_k$, 
$\hat c^{\alpha}$ 
$\hat f^{\alpha}$ and $\hat g$ in place of 
$a^{\alpha}_k$, $b^{\alpha}_k$, $c^{\alpha}$,   
$f^{\alpha}$ and $g$, respectively. We denote 
these solutions by $u$ and $\hat u$, respectively.

\begin{lemma}                                   \label{lemma 3.02.01.08}
Let Assumption \ref{assumption 1.01.01.09} 
hold. 
Let $g$ and $\hat g$ be bounded functions 
on $\mathbb R^d$.                           
Let $\e\in (0,Kh]$ be a constant 
and assume that for all $\alpha \in A$
  \begin{equation}                        \label{2.01.01.09}
    \begin{split}
      |b^\alpha_k - \hat{b}^\alpha_k| 
+ m^\alpha|f^\alpha - \hat{f^\alpha}| 
+ m^{\alpha}|c^\alpha - \hat{c}^\alpha| \leq K\e,          \\
      |a^\alpha_k - \hat{a}^\alpha_k| 
\leq K \e \sqrt{a^\alpha_k \land \hat{a}^\alpha_k} + K\e h,
    \end{split}
  \end{equation}
  on $\mt$. Then there is a constant $N^{\ast}$ depending 
on $K$ and $d_1$ such that for any constant 
$c_o\geq0$ satisfying \eqref{assumption-on-c_0} we have  
$$
    |u - \hat{u}| 
\leq  \e N_0e^{c_0 (T+\tau)}\big(N_1+\sup_{\mathbb R^d}(|g|+
|\hat g|
+\max_k|\delta_{h,\ell_k}g|
$$
\begin{equation*}
+\max_k|\delta_{h,\ell_k}\hat g|
+\varepsilon^{-1}|g-\hat g|)\big)
\end{equation*}
on $\bar{\mathcal M}_T$, where 
$N_0$ and $N_1$ are constants. The constant $N_0$ depends on $K$, $d_1$ and
$\rho$. The constant $N_1$ depends on $K$, $d_1$, $\rho$ and $\lambda$, 
provided $\lambda>0$, and it depends on 
$K$, $d_1$, $\rho$ and $T$ when $\lambda=0$. 
\end{lemma}

\begin{proof}
  We follow the idea of
  \cite{krylov:rate:lipschitz} to obtain 
  this lemma from the gradient estimate
  \eqref{eqn-grad-est-col}.
We consider $\mathbb R^d$ as a subspace 
$\mathbb R^d\times\{0\}$ of $\mathbb R^d\times\mathbb R$, 
and the vectors $\ell_k$ are identified with 
$(\ell_k,0)\in\mathbb R^{d+1}$ for $k=\pm1,\dots\pm d_1$.   
Let $(t,x) =
(t,x',x^{d+1})
\in [0,T]
  \times \R^d \times \R$. 
Let $\ell =(0,\ldots,0,1) \in \R^{d+1}$.  
Set $\ell_{\pm (d_1+1)}=\pm \ell$, 
$\delta_{h_r,\ell_r}=\delta_{h,\ell_k}$, 
for $r=k=\pm1\dots,\pm d_1$, 
$\delta_{h_r,\ell_r}
=\delta_{\varepsilon,\ell_{\pm(d_1+1)}}$ 
for $r=\pm(d_1+1)$, and
$$
\bar{\mathcal M}_T(\varepsilon)
= \bar{\mathcal M}_T \times \{0, \pm \e, \pm 2\e,
\ldots\},                                                   
\quad 
{\mathcal{M}}_T(\e) :=  \bar{\mathcal{M}}_T(\e) 
\cap([0,T)\times\mathbb R^d\times\mathbb R). 
$$
Let
  \begin{equation*}
    \tilde{a}^\alpha_k(t,x',x^{d+1}) = \left\{
      \begin{array}{lcl}
        a^\alpha_k(t,x') & \text{if} & x^{d+1} > 0,          \\
        \hat{a}^\alpha_k(t,x') & \text{if} & x^{d+1} \leq 0,
      \end{array} \right.
  \end{equation*}
  and define $\tilde{b}^\alpha_k,
  \tilde{c}^\alpha_k$, $\tilde{f}^\alpha$, 
$\tilde g$ 
  and $\tilde{u}$ similarly. Then
  $\tilde{u}$ satisfies
  \eqref{eq:fd:a}-\eqref{eq:fd:b} with 
$\mathcal M_T(\varepsilon)$, 
  $\tilde{a}^\alpha_k, \tilde{b}^\alpha_k,
  \tilde{c}^\alpha$, 
$\tilde{f}^\alpha$ and $\tilde g$
  in place of $Q$, $a^\alpha_k$, $b^\alpha_k$,
  $c^\alpha$, $f^{\alpha}$ and $g$, respectively. 
  To apply Corollary
  \ref{col-grad-est}
  to $\tilde{u}$ we need to check 
Assumption \ref{ass-for-thm-5.2} with 
$\mathcal M_T(\varepsilon)$, 
  $\tilde{a}^\alpha_k, \tilde{b}^\alpha_k,
  \tilde{c}^\alpha$ and  
$\tilde{f}^\alpha$
  in place of $Q$, $a^\alpha_k$, $b^\alpha_k$,
  $c^\alpha$ and $f^{\alpha}$, respectively. 
Clearly this assumption with 
$r= \pm 1, \ldots, \pm d_1$
holds by virtue of 
Assumption \ref{assumption 1.01.01.09}. 
Since 
$\delta_{\varepsilon,-\ell}
=-T_{\varepsilon,-\ell}\delta_{\varepsilon,\ell}$, 
we need only show that it holds also
for 
$r=(d+1)$.  To this end notice that 
\begin{equation*}
\fddelta_{\e,\ell}\tilde\psi(t,x) 
=\left\{
\begin{array}{lcl}
0  & \text{if} & x^{d+1} \neq 0,         \\
\e^{-1}(\psi(t,x') -\hat\psi(t,x'))  
& \text{if} & x^{d+1} = 0
\end{array}\right.
\end{equation*}
with $a^\alpha_k$, $b^{\alpha}$, $c^{\alpha}$ 
and $f^{\alpha}$ in place of $\psi$. 
Moreover, due to \eqref{2.01.01.09}
  \begin{equation*}
    \e^{-1}|a^\alpha_k(t,x') - \hat{a}^\alpha_k(t,x')| 
\leq K\sqrt{\tilde{a}^\alpha_k(t,x',0)} + Kh.
  \end{equation*}
Thus $|\fddelta_{\e,\ell}\tilde b|\leq K$, 
$m^{\alpha}|\fddelta_{\e,\ell}\tilde c|\leq K$,  
$m^{\alpha}|\fddelta_{\e,\ell}\tilde f|\leq K$ on 
$\mathcal M_T(\varepsilon)$, and  
$$
|\fddelta_{\e,\ell}\tilde a^{\alpha}_k|
\leq K\sqrt{\tilde a^{\alpha}_k}+Kh 
\quad 
\text{on $\mathcal M_T(\varepsilon)$}. 
$$
Hence we get the lemma by 
using Corollary \ref{col-grad-est}. 
\end{proof}

\begin{theorem}         
                              \label{theorem 3.17.3.9}
  Let Assumptions 
\ref{assumption 1.23.12.08} through 
\ref{assumption 2.23.12.08} and Assumption 
\ref{thm-est-disc-grad} hold for 
$a^{\alpha}_k$, $b^{\alpha}_k$ 
$c^{\alpha}$ and $f^{\alpha}$ and also 
for $\hat a^{\alpha}_k$, 
$\hat b^{\alpha}_k$,  
$\hat c^{\alpha}$ and $\hat f^{\alpha}$ in place of 
$a^{\alpha}_k$, $b^{\alpha}_k$ 
$c^{\alpha}$ and $f^{\alpha}$, respectively.  
Let $g$ and $\hat g$ be bounded functions on 
$\bar H_T$ such that for all $t\in[0,T]$, $x,y\in\mathbb R^d$
$$
|g(t,x)|+|\hat g(t,x)|\leq K, 
\quad
|g(t,x)-g(t,y)|+|\hat g(t,x)-\hat g(t,y)|\leq K|x-y|. 
$$
Set 
$$
    \e = \sup_{\mathcal{M}_T, A, k}
\left(|\sigma^\alpha_k - \hat{\sigma}^\alpha_k|
 + |b^\alpha_k - \hat{b}^\alpha_k| 
+ m^{\alpha}|c^\alpha - \hat{c}^\alpha| +
m^\alpha|f^\alpha - \hat{f}^\alpha|  
+|g-\hat g|\right), 
$$
where $\sigma_k^{\alpha}=\sqrt{a^{\alpha}_k}$, 
$\hat\sigma_k^{\alpha}=\sqrt{\hat a^{\alpha}_k}$. 
  Assume that 
$u$ and  $\hat{u}$ satisfy \eqref{eq:fd:a}-\eqref{eq:fd:b},  
$w$ and $\hat w$ satisfy  \eqref{1.17.02.09}-\eqref{2.17.02.09}
with $\mathcal{M}_T$ in place of $Q$, and  
$\sigma, b, c, f,g $ and
$\hat{\sigma}, \hat{b}, \hat{c}, \hat{f}, \hat{g}$, 
in place of $\sigma, b, c, f,g $,  
respectively. 
Then there is a constant $N^{\ast}$ depending 
on $K$, $d_1$ and $\rho$ such that 
for any constant $c_o\geq0$ satisfying 
\eqref{assumption-on-c_0} we have 
\begin{equation}                                        \label{eqn-4.5-ass}
    |u - \hat{u}| \leq Ne^{c_0( T+\tau)}\e, 
\quad
 |w - \hat{w}| \leq Ne^{c_0( T+\tau)}\e
\quad
\text{on $\bar{\mathcal M}_T$}, 
  \end{equation}
where $N$ is a constant depending on 
$K$, $d_1$, $\rho$ and $\lambda$, provided $\lambda>0$. 
If $\lambda=0$ then $N$ depends on $K$, $d_1$, $\rho$ and $T$. 
\end{theorem}

\begin{proof}
Consider first the case $\e \in (0,h]$. Then 
$$
|\sigma^\alpha_k -\hat{\sigma}^\alpha_k|\leq \varepsilon, 
\quad 
|\sigma^\alpha_k -\hat{\sigma}^\alpha_k|^2\leq \varepsilon h,  
$$
and using the identity
$$
    |a^2 - b^2| = (a+b)|a-b| = 2(a \land b)|a-b| + |a-b|^2,
$$
valid for any nonnegative numbers $a$ and $b$, we get
$$
    |a^\alpha_k - \hat{a}^\alpha_k| 
=2(|\sigma^\alpha_k| \land |\hat{\sigma}^\alpha_k|)
|\sigma^\alpha_k -\hat{\sigma}^\alpha_k| 
+ |\sigma^\alpha_k - \hat{\sigma}^\alpha_k|^2 
\leq2\e\sqrt{a^\alpha_k \land \hat{a}^\alpha_k} + \e h.
$$
  Hence by Lemma
  \ref{lemma 3.02.01.08},
  $|u-\hat{u}| \leq \e Ne^{c_0 (T+\tau)}$ on $\bmt$.
  Now consider the case $\e > h$. For
  $\theta \in [0,1]$, let $u^\theta$ be
  the solution of
  \begin{eqnarray*}
    \sup_{\alpha}m^\alpha(\fddeltatau u^\theta 
+ a^{\theta\alpha}_k \fdDeltak u^\theta + b^{\theta\alpha}_k \fddeltak
u^\theta  - c^{\theta\alpha}u^\theta + f^{\theta\alpha}) = 0 \on
\mathcal{M}_T\\
    g^\theta = u^\theta \on \{(T,x) \in \bar{\mathcal{M}}_T\},
  \end{eqnarray*}
  where 
  \begin{equation*}
    (\sigma^{\theta\alpha}_k, 
b^{\theta\alpha}_k, c^{\theta\alpha}, f^{\theta\alpha}, g^\theta) =
(1-\theta)(\sigma^\alpha_k, b^\alpha_k, c^\alpha, f^\alpha, g) +
\theta(\hat{\sigma}^\alpha_k, \hat{b}^\alpha_k, 
\hat{c}^\alpha, \hat{f}^\alpha, \hat{g}) 
\end{equation*}
  and $a^{\theta\alpha}_k = (1/2)|\sigma^{\theta\alpha}_k|^2$. 
For any $\theta_1, \theta_2 \in [0,1]$,
$$
    |\sigma^{\theta_1\alpha}_k - \sigma^{\theta_2 \alpha}_k| 
+ |b^{\theta_1\alpha}_k - b^{\theta_2 \alpha}_k| + |c^{\theta_1\alpha} -
c^{\theta_2 \alpha}|
$$
$$ 
+ m^\alpha|f^{\theta_1\alpha} - f^{\theta_2 \alpha}|
+|g^{\theta_1} - g^{\theta_2}|\!\! \leq\! |\theta_1 - \theta_2|\e. 
$$
  Hence if $\theta_1, \theta_2$ satisfy $|\theta_1 - \theta_2|\e 
\leq h$, then, thanks to the first part of the proof, with $u^{\theta_1}$ 
and $u^{\theta_2}$ playing the roles   
of $u$ and $\hat{u}$, respectively, 
  \begin{equation*}
    |u^{\theta_1} - u^{\theta_2}| \leq N |\theta_1 - \theta_2|
\varepsilon
e^{c_o(T+\tau)}.
  \end{equation*}
Set $\theta_i:=i/m$ for $i=0,1,\dots,m$ for an integer $m\geq1$ such that 
$\varepsilon/m\leq h$. 
Then 
\begin{equation*}                             \label{1.17.3.9}
|\hat u - u|\leq
 \sum_{i=0}^{m-1}|u^{\theta_{i+1}} - u^{\theta_i}|
\leq N \sum_{i=0}^{m-1}|\theta_{i+1} - \theta_i|\varepsilon
e^{c_oT^{\prime}}
= N\varepsilon e^{c_o(T+\tau)}, 
\end{equation*} 
that proves 
\eqref{eqn-4.5-ass} for $u$ and $\hat u$. Hence 
by using Remark \ref{remark 3.17.02.09} to rewrite 
equation \eqref{1.17.02.09}  
we get \eqref{eqn-4.5-ass} also for $w$ and $\hat w$.  
\end{proof}

\section{Some properties of the reward functions}             \label{section-analytic-ppties}

Let $A$ be a separable metric space.  
Let $\sigma=\sigma^{\alpha}(t,x)$ and 
$\beta=\beta^{\alpha}(t,x)$ be some 
Borel functions of 
$(\alpha,t,x)\in  A\times[0,\infty)\times\mathbb R^d$ 
with values in $\mathbb R^{d\times d^{\prime}}$ and $\mathbb R^d$, 
respectively.   
Let $\alpha=(\alpha_t)_{t\geq0}$ be a progressively 
measurable process with values in $A$, 
such that for every 
$s\in[0,T)$ and $x\in\mathbb R^d$ there is  
a solution $x_t=\{x_t^{\alpha,s,x}:t\in[0,T-s]\}$ of 
equation \eqref{eqn-sde}. 

Let $f=f^{\alpha}(t,x)$, 
$c=c^{\alpha}(t,x)\geq\lambda$ and $g=g(t,x)$ be Borel 
functions of 
$(\alpha,t,x)\in A\times[0,\infty)\times\mathbb R^d$ 
and of $(t,x)\in[0,\infty)\times\mathbb R^d$, respectively, where 
$\lambda\geq0$ is some constant. Set 
\begin{equation*}                             \label{6.15.01.09}
v^{\alpha}(s,x)
=E\int_0^{T-s}f^{\alpha_t}(s+t,x_t^{\alpha,s,x})
e^{-\varphi_{t}}\,dt
+Eg(T,x_{T-s}^{\alpha,s,x})e^{-\varphi_{T-s}}, 
\end{equation*}
\begin{equation}                               \label{1.18.3.9.}                         
w^{\alpha,\tau}(s,x)
=E\int_0^{\tau}f^{\alpha_t}(s+t,x_t^{\alpha,s,x})
e^{-\varphi_{t}}\,dt
+Eg(s+\tau,x_{\tau}^{\alpha,s,x})e^{-\varphi_{\tau}}, 
\end{equation}
$$
\varphi=\varphi_t^{\alpha,x}
=\int_0^tc^{\alpha_u}(s+u,x_u^{\alpha,s,x})\,du
$$
for $s\in[0,T]$, $x\in\mathbb R^d$,  
for the process $\alpha=(\alpha_t)$ and 
for a fixed stopping times $\tau$ with values in $[0,T-s]$.

\begin{lemma}                           \label{lemma 1.14.01.09} 
Assume that there exists a constant $K$  
such that 
$|g|\leq K$ on $\bar H_T$ 
and 
\begin{equation}                         \label{1.14.01.09}              
|f^{\alpha}(t,x)|\leq K(1+c^{\alpha}(t,x))
\end{equation}
for all $\alpha\in A$, $t\geq0$ and $x\in\mathbb R^d$. Then 
for $u:=v^{\alpha}, w^{\alpha,\tau}$ we have 
$$
|u|\leq K(2+N)\quad\text{on $\bar H_T$}, 
$$
where $N=(1-\exp(-\lambda T))/\lambda$ if $\lambda>0$, and $N=T$ 
if $\lambda=0$.  
\end{lemma}
\begin{proof}
Notice that 
$$
\int_0^{T-s}c^{\alpha}(s+t, x_{t})
e^{-\varphi_t}\,dt=1-e^{-\varphi_{T-s}}\leq1. 
$$
Hence 
$$
|u(s,x)|
\leq KE\int_0^{T-s}(1+c^{\alpha}(s+t, x_{t}))
e^{-\varphi_t}\,dt+K\leq K(2+N). 
$$
\end{proof}

\begin{assumptions}           \label{assumptions-for-properties-of-payoff-function}
  There exist a Borel function $m:A \to (0,1]$ 
and constants $\rho>0$,  $K\geq0$ and $L$
  such that for all $\alpha \in A$,
  $t\in[0,T]$ and  $x,y\in \R^d$ 
\begin{gather}
    m^\alpha(1+ c^\alpha(t,x)-\lambda) 
\geq \rho,\quad |m^\alpha f^\alpha(t,x)| 
\leq K,
                                                       \label{eqn-5.1}\\ 
    m^\alpha|f^\alpha(t,x) - f^\alpha(t,y)|\leq K|x-y| \label{3.14.01.09}\\
|c^\alpha(t,x) - c^\alpha(t,y)|\leq K|x-y|,           \label{eqn-5.2}\\
|g(t,x)-g(t,y)|\leq K|x-y|,                                   \label{4.14.01.09}\\ 
(x-y)(\beta^\alpha(t,x) - \beta^\alpha(t,y))
+ \tfrac{1}{2}|\sigma^\alpha(t,x) - \sigma^\alpha(t,y)|^2
\leq L|x-y|^2.                                         \label{eqn-5.3}
\end{gather}
 
\end{assumptions}
\begin{remark}
Notice that condition \eqref{eqn-5.1} 
implies condition \eqref{1.14.01.09} of 
Lemma \ref{lemma 1.14.01.09}, with $K/\rho$ 
in place of $K$ in \eqref{1.14.01.09}. 
Clearly, if $\beta^{\alpha}$ and $\sigma^{\alpha}$ are 
Lipschitz continuous in $x\in\mathbb R^d$, with 
Lipschitz constant $L/2$, independent of $\alpha\in A$ 
and $t\in[0,T]$, then the {\it monotonicity condition} 
\eqref{eqn-5.3} is satisfied. 
\end{remark}

\begin{lemma}                                              \label{lemma-lipschitz-continuity}
Let Assumption
  \ref{assumptions-for-properties-of-payoff-function}
hold. Assume 
\begin{equation}                   \label{condition 5.14.01.09}
|g|\leq K\quad \text{on $\bar H_T$}. 
\end{equation} 
Then for $u:=v^{\alpha}, w^{\alpha,\tau}$ we have 
  \begin{equation*}
    |u(s,x) -u(s,y)| 
\leq N|x-y|
\quad
\text{for all $s\in[0,T]$ and $x,y\in\mathbb R^d$}, 
  \end{equation*}
where $N$ is a constant depending only on $K$, $\rho$ and $T$. 
If $\lambda\geq |L|+2$, then $N$ depends only on $K$ and $\rho$. 
\end{lemma}

\begin{proof}
Clearly,  $|u(s,x) - u(s,y)|\leq \sum_{k=1}^4 I_k$
with  
\begin{equation*}
\begin{split}
& I_1=
\int_0^{T-s}|f^{\alpha_t}(s+t,x^{\alpha,s,x}_t)|
|e^{-\varphi_t^{\alpha,s,x}}-e^{-\varphi_t^{\alpha,s,y}}|dt,        \\
& I_2= E \int_0^{T-s} 
|f^{\alpha_t}(s+t,x^{\alpha,s,x}_t) -
f^{\alpha_t}(s+t,x^{\alpha,s,y}_t)|e^{-\varphi_t^{\alpha,s,y}} dt,          \\
&I_3=\sup_{\tau\in\ST(T-s)}E\{|g(s+\tau,x^{\alpha,s,x}_{\tau})|
|e^{-\varphi_{\tau}^{\alpha,s,x}}-e^{-\varphi_{\tau}^{\alpha,s,y}}|\},  \\
  &I_4 =\sup_{\tau\in\ST(T-s)}
E\{ |g(s+\tau,x^{\alpha,s,x}_{\tau}) - g(s+\tau,x^{\alpha,s,y}_{\tau})
|e^{-\varphi_{\tau}^{\alpha,s,y}}\}.
\end{split}
\end{equation*}
By \eqref{eqn-5.2} and \eqref{eqn-5.1}
\begin{equation*}
I_1 \leq E \int_0^{T-s}|f^{\alpha_t}(s+t,x_t^{\alpha,s,x})|
|\varphi_t^{\alpha,s,x} -
\varphi_t^{\alpha,s,y}|
e^{-\min(\varphi_t^{\alpha,s,x},\varphi_t^{\alpha,s,y})}dt
  \end{equation*}
$$
\leq  
K^2
E\int_0^{T-s}te^{-(\lambda-1)t}\tfrac{1}{m^{\alpha_t}}
\sup_{r\leq t}|x^{\alpha,s,x}_r-x^{\alpha,s,y}_r
|e^{-\int_0^t \tfrac{\rho}{m^{\alpha_u}}\,du} \,dt         \\
$$
$$
\leq 
N_1
\sup_{t\leq T-s}e^{-N_0t}
E\sup_{r\leq t}|x^{\alpha,s,x}_r-x^{\alpha,s,y}_r|, 
$$
for any constant $N_0\geq0$, 
where $N_1=K^2(e(\lambda-1-N_0))^{-1}\rho^{-1}$ 
when $\lambda\geq1+N_0$, and 
$N_1$ depends on $\rho$, $K$ $N_0$ 
and $T$ when $\lambda\in[0,1+N_0]$. 
 By \eqref{3.14.01.09} and \eqref{eqn-5.1} 
  \begin{equation*}
I_2\leq K
E\int_0^{T-s}(m^{\alpha_t})^{-1}
|x^{\alpha,s,x}_t -x^{\alpha,s,y}_t|
e^{-\varphi^{\alpha,s,y}_t} dt
  \end{equation*}
  \begin{equation*}
    \leq KE 
\left(\sup_{t\leq T}e^{-N_0t}
|x^{\alpha,s,x}_t-x^{\alpha,s,y}_t| 
\int_0^{T-s}e^{-(\lambda-1-N_0)t}
\tfrac{1}{m^{\alpha_t}}
e^{-\int_0^t\tfrac{\rho}{m^{\alpha_{u}}}\,du} dt 
\right) 
  \end{equation*}
  \begin{equation*}
    \leq N_2 E \sup_{t\leq T-s}
e^{-N_0t}|x^{\alpha,s,x}_t-x^{\alpha,s,y}_t|, 
  \end{equation*}
for every constant $N_0\geq0$, 
where $N_2=K/\rho$ when $\lambda\geq1+N_0$, 
and $N_2$ depends 
on $K$, $\rho$, $N_0$ and $T$ when $\lambda<1+N_0$. 
Due to conditions \eqref{condition 5.14.01.09}, 
\eqref{eqn-5.2}, $c^{\alpha}\geq\lambda$ and 
\eqref{4.14.01.09} we have   
\begin{equation*}
    I_3 \leq K\sup_{\tau\in\ST(T-s)}
E|e^{-\varphi_{\tau}^{\alpha,s,x}}
-e^{-\varphi_{\tau}^{\alpha,s,y}}|
\leq K E\sup_{t\leq T-s}e^{-\lambda t}
\int_0^t|x^{\alpha,s,x}_r-x^{\alpha,s,y}_r|\,dr
 \end{equation*}
$$
\leq K\int_0^{T-s}Ee^{-\lambda r}
|x^{\alpha,s,x}_r-x^{\alpha,s,y}_r|\,dr
\leq N_3\sup_{t\leq T-s}
Ee^{-N_0t}|x^{\alpha,s,x}_t-x^{\alpha,s,y}_t|
$$
for every constant $N_0\geq0$, where $N_3$ 
depends only on $K$ if 
$\lambda>N_0+1$ and $N_3$ depends on $K$, $N_0$ 
and $T$ if  
$\lambda\leq N_0+1$. 
Similarly, 
\begin{equation*}
 I_4 \leq K E \sup_{t\leq T-s}
e^{-\lambda t}|x^{\alpha,s,x}_t-x^{\alpha,s,y}_t|
\leq N_4 E \sup_{t\leq T-s}
e^{-N_0t}|x^{\alpha,s,x}_t-x^{\alpha,s,y}_t|
  \end{equation*}
for any $N_0\geq0$, where $N_4=K\exp((N_0-\lambda)T)$. 
  Consequently,
\begin{equation}                                    \label{1.14.02.09}
|u(s,x) - u(s,y)| 
\leq 
N E \sup_{t\leq T-s}e^{-N_0t}|x^{\alpha,s,x}_t-x^{\alpha,s,y}_t|   
\end{equation}
for every $N_0\geq0$. The constant 
 $N$ depends only on $K$ and $\rho$, 
if $\lambda\geq N_0+2$, and  it depends on $K$, $\rho$, 
$N_0$ and $T$ if $\lambda<2+N_0$. 
Using It\^o's formula and condition \eqref{eqn-5.3}, we have 
$$
e^{-2Lt}|x^{\alpha,s,x}_t-x^{\alpha,s,y}_t|^2\leq |x-y|^2+M_t
$$
almost surely for all $t\in[0,T-s]$, where 
$M$ is a local martingale. Thus 
$$
Ee^{-2L\tau}|x^{\alpha,s,x}_{\tau}-x^{\alpha,s,y}_{\tau}|^2\leq |x-y|^2
$$
for all stopping times $\tau\leq T-s$, that yields
$$
E\sup_{t\leq T-s}
e^{-Lt}|x^{\alpha,s,x}_{t}-x^{\alpha,s,y}_{t}|\leq 3|x-y|
$$
by virtue of Lemma 3.2 from \cite{GK2}. 
Combining this with estimate \eqref{1.14.02.09} we 
finish the proof of the lemma. 
\end{proof}

Assume that $A=\cup_{n=1}^{\infty}A_n$ 
for an increasing 
sequence of Borel sets $A_n$ of $A$ such that 
Assumptions \ref{assumption 1.15.01.09} 
and \ref{assumption 2.15.01.09} hold with $A_n$.  
Then the reward functions $v^{\alpha}$ 
and $w^{\alpha,\tau}$  
are well-defined on $\bar H_T$ for 
every $\alpha\in \A=\cup_{n=1}^{\infty}\A_n$, 
where $\A_n$ denotes the set of 
progressively measurable processes $\alpha=(\alpha_t)_{t\geq0}$ 
taking values in $A_n$. Thus   
we can define the optimal reward functions 
$$
v(s,x)=\sup_{\alpha\in\A}v^{\alpha}(s,x),\quad 
w(s,x)=\sup_{\alpha\in\A}
\sup_{\tau\in\ST(T-s)}w^{\alpha,\tau}(s,x) 
$$
for every $(s,x)\in[0,T]\times\mathbb R^d=\bar H_T$.
Recall the notation $H_T := [0,T) \times R^d$, 
$$
L^{\alpha}=\sigma^{\alpha}_{ik}\sigma^{\alpha}_{jk}
D_iD_j+\beta^{\alpha}_iD_i+c^{\alpha},   
$$
and let $C^{1,2}(\bar H_T)$ denote the set of functions 
$\psi=\psi(t,x)$ whose first derivative in $t$ and 
second order derivatives in $x$ are continuous functions 
on $\bar H_T$. 
The following lemma formulates an important 
property of smooth supersolutions and subsolutions 
to Bellman equations. 

\begin{lemma}                                      \label{lemma  1.11.02.09}
Let Assumptions \ref{assumption 1.15.01.09} 
and 
\ref{assumption 2.15.01.09} hold.  Assume that 
$\sigma$, $\beta$  
are continuous in $\alpha\in A$. 
Assume, moreover,  that  
$f$ and $c$ are continuous 
in $(\alpha,x)$ and  
are continuous in $x$, uniformly in $\alpha\in A$, 
for each $t\in[0,T]$.     
Let $S\in(0,T]$ and 
$\psi \in C^{1,2}(\bar H_S)$ 
such that for some constants 
$K$ and $q\geq0$ 
\begin{equation}                                 \label{condition 8.15.01.09}
|\psi(t,x)|\leq K(1+|x|^q)
\quad
\text{for all $(t,x)\in H_S$}.  
\end{equation}
Let  
$Q$ be a domain  contained in $H_S$. Denote   
its boundary by $\partial Q$. 
Then the following statements hold:  
\begin{enumerate}
\item[(i)] Let 
  \begin{equation}                               \label{condition 2.03.02.09}
    \tfrac{\partial}{\partial  t}\psi 
+ L^{\alpha} \psi + f^{\alpha} \leq 0 
\quad
\text{on $Q$, for all $\alpha\in A$}.
  \end{equation} 
 Then 
\begin{equation}                                     \label{eq:supersoln}
v \leq \psi + \sup_{\partial Q}
[v - \psi]_+ 
\quad \text{on $\bar Q$}.  
\end{equation} 
In addition to \eqref{condition 2.03.02.09} let $g\leq \psi$ on $Q$. 
Then \eqref{eq:supersoln} holds also for $w$ in place of $v$.  
\item[(ii)]
 Let 
  \begin{equation}                               \label{condition 3.11.02.09}
    \tfrac{\partial}{\partial  t}\psi 
+ L^{\alpha} \psi + f^{\alpha} \geq 0 
\quad
\text{on $Q$, for some $\alpha\in A$}.
  \end{equation} 
 Then 
\begin{equation}                                     \label{5.11.02.9}
v \geq \psi - \sup_{\partial Q}
[v - \psi]_- 
\quad\text{and\quad $w \geq \psi - \sup_{\partial Q}
[w - \psi]_-$ on\quad $\bar Q$} .  
\end{equation} 
\end{enumerate}

\end{lemma}

\begin{proof}
This lemma follows from Lemma 6.1.2 
and 
Theorem 6.1.5 from \cite{krylov:controlled}. 
For the convenience of the reader we give 
a more detailed proof here. 
Set
$v_n=\sup_{\alpha\in\A_n}v^{\alpha}$  for integers $n\geq1$.    
Then by Theorem 3.1.5 in 
\cite{krylov:controlled}, 
the polinomial growth condition  
\eqref{condition 8.15.01.09} holds for  
$v_n$ in place of $\psi$, with some constants 
$K$ and $q$ depending on $n$,   
and $v_n$ is continuous on $\bar H_T$. 
Set 
$$
\tau_{Q}=\inf\{t\geq0:(s+t,x_t)\notin Q\},
\quad \tau^R_Q=\inf\{t\geq0:|x_t|\geq R\}\wedge \tau_Q 
$$
for $R>0$. By Bellman's principle (Theorem 2.3.6 from 
\cite{krylov:controlled}),  for $(s,x)\in Q$, 
integer $n\geq1$, stopping time $\tau=\tau_Q^R$, for 
any $\e > 0$ there is a strategy 
$(\alpha_t) \in \A_n$ such that
\begin{equation}                    \label{1.18.01.09}
v_n(s,x) \leq \varepsilon + I^{(\alpha)}_n(s,x), 
\end{equation}                      
\begin{equation}                    \label{12.11.02.09}                           
I^{(\alpha)}_n(s,x):=E_{s,x}^{\alpha} 
\left(\int_0^{\tau} 
f^{\alpha_t}(s+t,x_t)e^{-\varphi_t} dt 
+ v_n(s+\tau,x_{\tau})e^{-\varphi_{\tau}}\right),  
\end{equation}
where, as before, to ease notation we use $x_t$ in place of 
$x^{\alpha,s,x}_t$. 
Using condition \eqref{condition 2.03.02.09} 
and applying 
It\^o's formula to
  $\psi(s+t,x_t)e^{-\varphi_t}$ 
we have   
$$
I^{(\alpha)}_n(s,x)\leq- E_{s,x}^{\alpha}\int_0^{\tau}
e^{-\varphi_t}
\left[\tfrac{\partial}{\partial t} + L^{\alpha_t}
\right]
\psi(s+t,x_t) dt  + E_{s,x}^{\alpha}
v_n(s+\tau,x_{\tau})e^{-\varphi_{\tau}} 
$$
$$
=\psi(s,x)
+E_{s,x}^{\alpha}
\{(v_n(s+\tau,x_{\tau})-\psi(s+\tau,x_{\tau}))
e^{-\varphi_{\tau}}\}. 
$$
Letting here $R\to\infty$ we get 
$$
I^{(\alpha)}_n(s,x)\leq 
\psi(s,x)
+E_{s,x}^{\alpha}
\{(v_n(s+\tau_Q,x_{\tau_Q})-\psi(s+\tau_Q,x_{\tau_Q}))
e^{-\varphi_{\tau_Q}}\}.     
$$
Thus from \eqref{1.18.01.09} we have 
$$
v_n(s,x) \leq \e + \psi(s,x) +
 \sup_{\partial Q} 
\left[v_n- \psi\right]_+ 
\leq \e + \psi(s,x) +
 \sup_{\partial Q} \left[v - \psi\right]_+.
$$
Letting here $n\to\infty$ and  $\e \to 0$ we get 
\eqref{eq:supersoln}. 
Hence \eqref{eq:supersoln} is valid also 
for $w$ in place of $v$, since 
$w=\sup_{\gamma\in\bar\A}v^{\gamma}$ 
by virtue of Theorem \ref{thm-opt-stop-and-cont-equiv}, 
Assumptions 
\ref{assumption 1.15.01.09}-\ref{assumption 2.15.01.09} 
remain valid with $\bar A_n$ and
$\bar A$  in place of $A_n$ and $A$, and due to  
\eqref{condition 2.03.02.09} and $\psi\geq g$ on $Q$,  
$$
\tfrac{\partial}{\partial t}\psi+L^{\gamma}\psi+f
=\tfrac{\partial}{\partial t}\psi+L^{\alpha}\psi+f +r(g-\psi)\leq 0
\quad\text{on $Q$}
$$
for every $\gamma=(\alpha,r)\in \bar A$. 
To prove (ii) let $\alpha\in A$ such that 
\eqref{condition 3.11.02.09} holds. Then $\alpha\in A_n$ 
for some $n\geq1$, the constant strategy $\alpha_t=\alpha$ 
belongs to $\A_n$, and by Bellman's principle
$$
v_n(s,x)\geq I^{(\alpha)}_n(s,x) 
$$
with this strategy $\alpha$, where $I^{(\alpha)}_n$ 
is defined by \eqref{12.11.02.09}. Hence by an obvious 
modification of the proof of part (i) we get the first 
inequality in \eqref{5.11.02.9}, and that 
yields the second inequality 
by virtue of Theorem \ref{thm-opt-stop-and-cont-equiv}, since 
clearly 
$\tfrac{\partial}{\partial t}\psi+L^{\gamma}\psi+f^{\gamma}\geq0$ 
on $Q$
for $\gamma=\alpha\in A\subset\bar A$.  
\end{proof}

Next we want to study the regularity of $v$ and 
$w$ in $t\in[0,T]$. 
The following simple example shows 
that Assumption \ref{assumptions-for-properties-of-payoff-function}
does not ensure the continuity of $v$ at $t=T$,  
even if 
$\sigma^{\alpha}$ and $b^{\alpha}$ are as regular as we wish.

\begin{example}
  \label{ex:hc}
  Let $A=[0,\infty)$, $f^{\alpha}(t,x)=\alpha$, $g(x)=0$, 
$c^{\alpha}(t,x)=\alpha$ for $\alpha\in A$. Then 
Assumption \ref{assumptions-for-properties-of-payoff-function} 
holds with 
$m^{\alpha}=(1+\alpha)^{-1}$ and $\sigma^{\alpha}=0$, $b^{\alpha}=0$,   
and for $t\in[0,T]$, $x\in\mathbb R^d$  
  \begin{equation*}
    v(t,x) = \sup_{\alpha \in \A} 
E\int_0^{T-t} \alpha_se^{-\int_0^s \alpha_u du}ds=
  \sup_{\alpha \in \A} E\left(1 - e^{-\int_0^{T-t} \alpha_u du}\right) 
= \one_{t<T},
  \end{equation*}
  which is not continuous at $T$.
\end{example}

\section{H\"older continuity in time}                             \label{section-hc}

Let $\sigma=\sigma^{\alpha}(t,x)$, $\beta=\beta^{\alpha}(t,x)$, 
$f=f^{\alpha}(t,x)$ and $c=c^{\alpha}(t,x)$  
be Borel functions of 
$(\alpha,t,x)\in A\times \mathbb R_+\times \mathbb R^d$,  
taking values in  $\mathbb R^{d\times d^{\prime}}$, $\mathbb R^d$, 
$\mathbb R$ and $\mathbb R_+$, respectively, 
such that $c\geq\lambda$ for a constant $\lambda\geq0$. 
Let $g$ be a Borel function on $\mathbb R_+\times\mathbb R^d$ 
with values in $\mathbb R$. 

We make the following assumption. 

\begin{assumptions}                         \label{assumption 5.03.02.09}
There is a constant $K$ such that for 
$\psi=\sigma^{\alpha}, \beta^{\alpha}, f^{\alpha}, c^{\alpha},g$  
for all $\alpha\in A$ we have 
$$
|\psi(t,x)-\psi(t,y)|\leq K|x-y|, 
\quad 
|\psi(t,x)|\leq K
$$
for all $t\in[0,T]$ and $x, y\in\mathbb R^d$.
\end{assumptions}

Obviously Assumption \ref{assumption 5.03.02.09} 
implies Assumptions \ref{assumption 1.15.01.09} and 
\ref{assumption 2.15.01.09}, the reward functions 
$v^{\alpha}$, $w^{\alpha,\tau}$, 
$v$ and $w$ are well-defined 
by \eqref{5.17.01.09}, \eqref{1.03.02.09} 
and  \eqref{eq:stop-payofffn}. Moreover, 
Assumption 
\ref{assumptions-for-properties-of-payoff-function} 
holds with $m^{\alpha}=\rho=1$ and $L=2K$. 
Thus by 
Lemma \ref{lemma-lipschitz-continuity} 
there is a constant $C$  
such that for $u:=v,w$
\begin{equation}                               \label{6.03.02.09}
|u(t,x)-u(t,y)|\leq C|x-y|
\quad \text{for all $t\in[0,T]$ and $x,y\in\mathbb R^d$}. 
\end{equation}
If $\lambda\geq K+2$, then $C$ depends only on $K$, otherwise 
it depends on $K$ and $T$. 
Using results from \cite{krylov:controlled} and 
\cite{krylov:rate:lipschitz:published}
one can prove the following lemma on the H\"older continuity 
of $v$ and $w$ in $t$. 

\begin{lemma}                                 \label{lemma 1.10.02.09}
Let Assumption  \ref{assumption 5.03.02.09} hold. 
Assume that 
$\sigma$, $\beta$  
are continuous in $\alpha\in A$. 
Assume, moreover,  that  
$f$ and $c$ are continuous 
in $(\alpha,x)$ and  
are continuous in $x$, uniformly in $\alpha\in A$, 
for each $t\in[0,T]$. Then for $x_0\in\mathbb R^d$ 
and $0\leq t_0\leq s_0\leq T$ such that $|s_0-t_0|\leq 1$, 
we have     
\begin{equation}                               \label{2.10.02.09}
|v(t_0,x_0)-v(s_0,x_0)|\leq N(\nu_1+1)|s_0 - t_0|^{1/2},  
\end{equation}
\begin{equation}                               \label{2.22.02.09}
|w(t_0,x_0)-w(s_0,x_0)|\leq N(\nu_2+1)|s_0 - t_0|^{1/2}
+\mu|s_0 - t_0|^{1/2},  
\end{equation} 
where $N$ is a constant depending only on $K$,  
and    
$$
\nu_1:=\sup_{y\in\mathbb R^d\setminus\{x_0\}}
\tfrac{|v(s_0,x_0)-v(s_0,y)|}{|x_0-y|}, 
\quad 
\nu_2:=\sup_{y\in\mathbb R^d\setminus\{x_0\}}
\tfrac{|w(s_0,x_0)-w(s_0,y)|}{|x_0-y|}, 
$$
\begin{equation}                                  \label{1.27.2.9}
\mu:=\sup_{y\in\mathbb R^d}\sup_{0\leq t<s_0}
\tfrac{|g(t,y)-g(s_0,y)|}{|t-s_0|^{1/2}}. 
\end{equation}
\end{lemma}
\begin{proof}
We may assume $\nu_1<\infty$, $\nu_2<\infty$, $\mu<\infty$ 
and $0\leq t_0<s_0$. Moreover, by
shifting  the origin we may assume $t_0 = 0$ and hence $s_0\leq 1$. 
To prove \eqref{2.10.02.09}
define for a constant $\gamma>0$ the function 
$$                                       
\psi(t,x) = \gamma\nu_1 [\xi(t)|x-x_0|^2 +\kappa_1(s_0 - t)] 
+ \kappa_2(s_0 -t) 
$$
\begin{equation}                                       \label{2.14.02.09} 
+ \nu_1\gamma^{-1} + v(s_0,x_0),
\quad\text{for $(t,x)\in \bar H_{s_0}$}, 
  \end{equation}
where $\xi(t)=\exp(s_0-t)$ and 
$\kappa_1>0$, $\kappa_2>0$ are 
some constants to be chosen later. 
By simple calculations for any $\alpha\in A$ 
$$
L^{\alpha}\psi(t,x)
=\nu_1\gamma \xi(t)[2\sigma^{\alpha}_{ik}\sigma^{\alpha}_{ik}(t,x)
+2\beta^{\alpha}_{i}(t,x)(x_i-x_{0i})]
$$
$$
 -c^{\alpha}(t,x)\psi(t,x)
\leq N_1\gamma\nu_1 (1+|x-x_0|)+N_2, 
$$
for $(t,x)\in\bar H_{s_0}$, where $N_1$ and $N_2$ 
are constants depending 
only on $K$. 
Hence, choosing $\kappa_2\geq K+N_2$,  
we have   
$$
\tfrac{\partial}{\partial t}\psi(t,x)
+L^{\alpha}\psi(t,x)+f^{\alpha}(t,x)
$$
\begin{equation}                                     \label{1.24.1.9}
\leq \nu_1\gamma[N_1(1+|x-x_0|)-|x-x_0|^2-\kappa_1], 
\end{equation}
where the right-hand side is negative for all 
$x$ if $\kappa_1$ is sufficiently large, depending only 
on $N_1$. 
Notice that for all $x\in\mathbb R^d$ 
$$
\psi(s_0,x)=\nu_1(\gamma |x-x_0|^2
+\gamma^{-1})+v(s_0,x_0)
\geq \nu_1|x-x_0|+v(s_0,x_0)\geq v(s_0,x). 
$$
Thus applying part (i) of Lemma \ref{lemma  1.11.02.09} 
with $S:=s_0$ 
and $Q:=H_{s_0}$ we obtain 
$$
v(t,x_0)\leq \nu_1[\gamma\kappa_1(s_0-t)+\gamma^{-1}]+\kappa_2(s_0-t)+v(s_0,x_0)
$$
for all $t\in[0,s_0]$ and constants $\gamma>0$. 
For $t=0$ we choose $\gamma=(\kappa_1 s_0)^{-1/2}$ to get
$$
v(0,x_0)\leq 2\nu_1\kappa^{1/2}_1s_0^{1/2}+\kappa_2s_0+v(s_0,x_0), 
$$
that yields 
\begin{equation}                                    \label{2.24.2.9}
v(0,x_0)-v(s_0,x_0)\leq N(\nu_1+1)s_0^{1/2}
\end{equation}
with $N=\max(2\kappa^{1/2}_1,\kappa_2)$.
To get the corresponding estimate 
for $w$, instead of \eqref{2.14.02.09} define 
$\psi$ by 
$$                                       
\psi(t,x) = \gamma \nu_2 [\xi(t)|x-x_0|^2 +\kappa_1(s_0 - t)] 
+ \kappa_2(s_0 -t) 
$$
\begin{equation}                                       \label{5.14.02.09} 
+ \nu_2\gamma^{-1} +\mu s_0^{1/2}+w(s_0,x_0). 
  \end{equation}
Then just like before we see that 
for sufficiently large constants $\kappa_1$ 
and  $\kappa_2$, depending only on $K$, 
the left-hand side of \eqref{1.24.1.9} remains 
negative for all $(t,x)\in\bar H_{s_0}$, 
and that 
$$
\psi(t,x)\geq\psi(s_0,x)\geq w(s_0,x) +\mu s_0^{1/2}\geq g(s_0,x)+\mu s_0^{1/2}
\geq g(t,x)
$$
for all $t\in [0,s_0]$ and $x\in\mathbb R^d$.    
Hence 
by part (i) of Lemma \ref{lemma  1.11.02.09}  
$$
w(0,x_0)\leq \nu_2[\gamma\kappa_1s_0+\gamma^{-1}]+\kappa_2s_0
+\mu s_0^{1/2}+w(s_0,x_0)
$$
for any $\gamma>0$, 
that yields  
\begin{equation*}                               \label{4.24.2.9}
w(0,x_0)\leq N(\nu_2+1)s_0^{1/2}+w(s_0,x_0).  
\end{equation*}
Now we prove this inequality with
$w(0,x)$ and $w(s_0,x)$ interchanged, together 
with inequality 
\eqref{2.24.2.9} with $v(0,x)$ and $v(s_0,x)$  
interchanged. To this end 
set 
$$
\psi(t,x)=-\gamma \nu [\xi(t)|x-x_0|^2 +\kappa_1(s_0 - t)] 
- \kappa_2(s_0 -t)  - C\gamma^{-1} + u(s_0,x_0). 
$$
with  $u:=\text{$v$ and $w$}$, and $\nu:=\text{$\nu_1$ and $\nu_2$}$, respectively.  
Notice that for large $\kappa_2$, depending only on $K$, 
we have 
$$
\tfrac{\partial}{\partial t}\psi(t,x)
+L^{\alpha}\psi(t,x)+f^{\alpha}(t,x)
\geq -\nu\gamma[N_1(1+|x-x_0|)-|x-x_0|^2-\kappa_1], 
$$
with a constant $N_1$ depending
on
$K$,  
where the right-hand side is positive for all 
$x$ if $\kappa_1$ is sufficiently large, depending only 
on $N_1$. 
Furthermore, 
$$
\psi(s_0,x)=-\nu(\gamma |x-x_0|^2+\gamma^{-1})+u(s_0,x_0)
$$
$$
\leq -\nu|x-x_0|+u(s_0,x_0)  
\leq u(s_0,x). 
$$
Hence by virtue of part (ii) of 
Lemma \ref{lemma  1.11.02.09} 
we get 
$$
u(t,x_0)
\geq -\nu[\gamma\kappa_1(s_0-t)+\gamma^{-1}]-\kappa_2(s_0-t)+u(s_0,x_0)
$$
for all $t\in[0,s_0]$ and constant $\gamma>0$. Choosing here 
$t=t_0=0$ and $\gamma=(\kappa_1 s_0)^{-1/2}$ we get 
$$
u(0,x_0)\geq -2\nu\kappa^{1/2}_1-\kappa_2s_0^{1/2}+u(s_0,x_0)
\geq -N(\nu+1)s_0^{1/2}+u(s_0,x_0)
$$
with $N:=\max(2\kappa^{1/2}_1,\kappa_2)$, 
that completes the proof of the lemma. 
\end{proof}
\begin{theorem}                                        \label{theorem 1.26.2.9}
Let Assumption \ref{assumption 5.03.02.09} 
hold. Assume that
$\sigma^{\alpha}$, 
$\beta^{\alpha}$, $f^{\alpha}$ are continuous in 
$\alpha\in A$ and that 
$$
|g(s,x)-g(t,x)|\leq K|t-s|^{1/2}
\quad\text{for all $t,s\in[0,T]$ and $x\in\mathbb R^d$}. 
$$
Then there is a constant $N$ such that for 
$u:=v,w$ we have 
$$
|u(s,x)-u(t,x)|\leq N|t-s|^{1/2}
\quad\text{for all $t,s\in[0,T]$ and $x\in\mathbb R^d$}. 
$$
The constant $N$ depends on $K$ and $T$. Moreover, 
there is a constant $\lambda_0$, depending on $K$,  
such that if $\lambda\geq\lambda_0$, 
then $N$ depends only on $K$. 
\end{theorem}
\begin{proof}
We get this theorem immediately from the previous lemma 
by taking into account Lemmas \ref{lemma 1.14.01.09} 
and \ref{lemma-lipschitz-continuity}. 
\end{proof}

Now we formulate the corresponding results for 
the solutions $v=v_{\tau,h}$ and $w=w_{\tau,h}$ 
of the finite difference schemes 
\eqref{2.30.12.08}-\eqref{3.30.12.08} and 
\eqref{1.18.02.09}-\eqref{2.18.02.09}, respectively, 
when $m^{\alpha}=1$ for all $\alpha\in A$.  
The following lemma is proved 
in \cite{krylov:rate:lipschitz:published}
for $u=v_{\tau,h}$.

\begin{lemma}                                                \label{lem:hold:tm:ct:vth}
 Let $\tau,h\leq K$. 
Let Assumption \ref{assumption 1.23.12.08}
hold and assume that for $\psi:=a_k^{\alpha}$, $b^{\alpha}_k$, 
$c^{\alpha}$,   
$f^{\alpha}$ and $g$ for every $k=\pm1,\dots\pm d_1$ 
and $\alpha\in A$ 
we have $|\psi|\leq K$ on $\bar H_T$. Let 
$(t_0,x_0)\in\bar H_T$ and $s_0\in[t_0,T]$ such that 
$s_0-t_0\leq 1$ and $(s_0-t_0)/\tau$ is an integer. 
Then \eqref{2.10.02.09} and \eqref{2.22.02.09} 
hold with  $v_{\tau,h}$ and $w_{\tau,h}$
in place of $v$ and $w$, respectively, 
where the constants $\nu_1$, 
$\nu_2$ and $\mu$ are defined  by \eqref{1.27.2.9}  
with  $v_{\tau,h}$ and $w_{\tau,h}$, 
in place of $v$ and $w$, respectively, and the constant 
$N$ depends on $K$ and $d_1$. 
\end{lemma}

\begin{proof}
We may assume that $s_0 > 0$ and also, by shifting 
the origin, that $t_0 = 0$, $x_0=0$ and hence that 
$s_0\in(0,1)$ is an integer multiple 
of $\tau$. Now we can prove the required estimates 
in the same way as Lemma \ref{lemma 1.10.02.09} 
is proved. We need only use 
Corollary \ref{corollary 1.25.2.9} with $T:=s_0$ and $Q:=\mathcal M_{s_0}$ 
in place of Lemma \ref{lemma  1.11.02.09}. 
\end{proof}

\begin{theorem}                                  \label{cor:hold:cont:in:t} 
Let $\tau,h\leq K$. 
Let Assumption \ref{assumption 1.23.12.08}
hold and assume that for $\psi:= 
\sqrt{a^{\alpha}_k}$, $b^{\alpha}_k$, 
$c^{\alpha}$,   
$f^{\alpha}$ and $g$,  
for every $k=\pm1,\dots\pm d_1$ 
and $\alpha\in A$ we have 
$$
|\psi(t,x)-\psi(s,y)|\leq K(|x-y|+|s-t|^{1/2}), 
\quad 
|\psi(t,x)|\leq K
$$
for all $s,t\in[0,T]$ and $x,y\in\mathbb R^d$. 
Then for $u:=v_{\tau,h}$,  $w_{\tau,h}$ we
have 
\begin{equation}                    \label{2.17.3.9}
    |u(t,x) - u(s,x)| 
\leq N\left(|t-s|^{1/2} + \tau^{1/2}\right)
  \end{equation}
for all $x\in\mathbb R^d$ and $s,t\in[0,T]$, 
where $N$ is a constant depending
only  on $K$, $d_1$ and $T$. There is a constant $\lambda_0\geq0$, 
depending only on $K$ and $d_1$, 
such that if $\lambda\geq\lambda_0$ then $N$ depends only 
on $K$ and $d_1$. 
\end{theorem}

\begin{proof} For $u=v_{\tau,h}$ estimate \eqref{2.17.3.9}
is proved in \cite{krylov:rate:lipschitz:published} 
(see Lemma 6.2 there). We get \eqref{2.17.3.9}
for $u=w_{\tau,h}$ similarly, noticing that 
Assumptions \ref{ass-new-for-3} and \ref{assumption 2.23.12.08} 
are obviously satisfied with $m^{\alpha}=1$ and 
$\rho=1$, and   
by using Lemma \ref{lem:hold:tm:ct:vth},  
Theorems \ref{theorem 4.17.3.9}, \ref{theorem 3.17.3.9} and  
Corollary \ref{corollary-soln-to-disc-bellman-pde-is-bdd}.  
\end{proof}

\section{Shaking and  Smoothing}                       
                                               \label{section-shaking}

The {\it method of shaking} 
is introduced in 
\cite{krylov:rate:variable}.  
Following \cite{krylov:rate:lipschitz:published} we adapt
it to optimal stopping of controlled  
diffusion processes and to the corresponding finite 
difference schemes. 

For $\varepsilon\in\mathbb R$ we set 
$$
A^{\varepsilon}=A\times[-\varepsilon^2, 0]
\times \{x\in\mathbb R^d:|x|\leq\varepsilon\},  
\quad
\bar A^{\varepsilon}=A^{\varepsilon}\times[0,\infty),   
$$
and identify $\alpha\in A$ with 
$(\alpha,0,0)\in A^{\varepsilon}$ and 
$(\alpha,\eta,\xi)\in A^{\varepsilon}$ 
with 
$(\alpha,\eta,\xi,0)\in \bar A^{\varepsilon}$. 
Thus $A\subset A^{\varepsilon}\subset \bar A^{\varepsilon}$. 

First we shake optimal stopping and control problems.  
Let
$\sigma=\sigma^{\alpha}(t,x)$,
$\beta=\beta^{\alpha}(t,x)$, 
$f=f^{\alpha}(t,x)$ and $c=c^{\alpha}(t,x)$  
be Borel functions of 
$(\alpha,t,x)\in A\times \mathbb R\times \mathbb R^d$,  
taking values in  $\mathbb R^{d\times d^{\prime}}$, $\mathbb R^d$, 
$\mathbb R$ and $\mathbb R_+$, respectively, 
such that $c\geq\lambda$ for a constant $\lambda\geq0$. 
Let $g$ be a Borel function on $\mathbb R\times\mathbb R^d$ 
with values in $\mathbb R$. 

We make the following assumption. 

\begin{assumptions}                         \label{assumption 5.17.3.9}
There is a constant $K$ such that for 
$\psi=\sigma^{\alpha}, \beta^{\alpha}, 
f^{\alpha}, g$, $c^{\alpha}-\lambda$,  
for all $\alpha\in A$ we have 
$$
|\psi(t,x)-\psi(s,y)|\leq K(|x-y|+|s-t|^{1/2}), 
\quad 
|\psi(t,x)|\leq K
$$
for all $s,t\in\mathbb R$ and $x, y\in\mathbb R^d$. 
\end{assumptions}

For $\gamma=(\alpha,\eta,\xi,r)\in \bar A^{\varepsilon}$ 
we set 
\begin{equation*}                                 \label{12.19.3.9}
\sigma^{\gamma}(t,x)=\sigma^{\alpha}(t+\eta, x+\xi), 
\quad 
\beta^{\gamma}(t,x)=\beta^{\alpha}(t+\eta, x+\xi), 
\end{equation*}
$$
c^{\gamma}(t,x)=c^{\alpha}(t+\eta, x+\xi)+r,
$$ 
\begin{equation}                                 \label{13.19.3.9}
f^{\gamma}(t,x)=f^{\alpha}(t+\eta, x+\xi)
+rg^{\varepsilon}(t,x), 
\end{equation}
for all $(t,x)\in\mathbb R\times\mathbb R^d$, where 
\begin{equation}                                    \label{10.19.3.9}
g^{\varepsilon}(t,x):=\sup_{\eta\in[-\varepsilon^2,0]}
\sup_{\xi\in\mathbb R^d,|\xi|\leq\varepsilon}|g(t+\eta,x+\xi)|.
\end{equation}
Let $\A^{\varepsilon}$ be the set of $A^{\varepsilon}$-valued 
progressively  measurable processes $(\gamma_t)_{t\geq0}$. 
Set 
$\bar\A^{\varepsilon}=\cup_{n=1}^{\infty}\bar\A^{\varepsilon}_n$, 
where $\bar\A^{\varepsilon}_n$ is the set of 
progressively measurable processes $(\gamma_t)_{t\geq0}$ 
with values in $\bar A^{\varepsilon}_n=A^{\varepsilon}\times[0,n]$.

Shaking the optimal reward $w$ given by 
\eqref{eq:stop-payofffn} means that we consider 
$\tilde w=\tilde w^{\varepsilon}(s,x)$ defined by 
\begin{equation*}                              
w^{\varepsilon}(s,x)=\sup_{\gamma\in\mathfrak{A}^{\varepsilon}}
\sup_{\tau\in\ST(T-s)}
w^{\gamma,\tau}(s,x),  
\end{equation*}
where 
\begin{equation*}                               
w^{\gamma,\tau}(s,x)=E^{\gamma}_{s,x}
\left[\int_0^{\tau}f^{\gamma_t}(s+t,x_t)
e^{-\varphi_t}dt + g^{\varepsilon}(s+\tau,x_{\tau})e^{-\varphi_{\tau}}
\right],
\end{equation*}
$$                                         
\varphi_t = \varphi_t^{\gamma,s,x}
= \int_0^t c^{\gamma_{r}}(s+r,x_{r}^{\gamma,s,x})\,dr.  
$$
Notice that if Assumption \ref{assumption 5.17.3.9} holds, 
then by virtue of 
Theorem \ref{thm-opt-stop-and-cont-equiv}
$$
w^{\varepsilon}=\sup_{\gamma\in\bar\A^{\varepsilon}}
v^{\gamma}=\lim_{n\to\infty}w_n^{\varepsilon}, 
$$
where
\begin{equation*}                              
v^{\gamma}(s,x)=E^{\gamma}_{s,x}
\left[\int_0^{T-s}f^{\gamma_t}(s+t,x_t)
e^{-\varphi_t}dt + g^{\varepsilon}(T,x_{T-s})e^{-\varphi_{T-s}}
\right],  
\end{equation*}
\begin{equation}                                \label{10.21.3.9}
w_n^{\varepsilon}
=\sup_{\gamma\in\bar\A^{\varepsilon}_n}v^{\gamma}. 
\end{equation}

\begin{lemma}                                      \label{lemma 1.18.3.9}
Let Assumption \ref {assumption 5.17.3.9} hold. 
Then there is a constant $N$ such that 
\begin{equation}                                     \label{7.18.3.9}
|w^{\varepsilon}-w|\leq N\varepsilon
\quad
\text{on $\bar H_T$}.  
\end{equation}
In addition to Assumption \ref {assumption 5.17.3.9} 
assume that $\sigma^{\alpha}$, 
$\beta^{\alpha}$ and $f^{\alpha}$ 
are continuous in $\alpha\in A$. 
Then there is a constant $N$ 
such that  
\begin{equation}                                       \label{8.18.3.9}
|w^{\varepsilon}(t,x)-w^{\varepsilon}(s,y)|
\leq N(|x-y|+|t-s|^{1/2})
\end{equation}
for $s,t\in[0,T]$, $x,y\in\mathbb R^d$.  
The constant $N$ in the above 
estimates depends only on $K$ and $T$. 
Moreover, there is a constant $\lambda_0$, 
depending only on $K$ such that $N$ 
is independent of $T$ 
if $\lambda\geq\lambda_0$. 
\end{lemma}

\begin{proof} 
Applying Lemma \ref{lemma-lipschitz-continuity} and Theorem 
\ref{theorem 1.26.2.9} we immediately 
get estimate \eqref{8.18.3.9}. 
By using the inequality 
$$
|a_1e^{-b_1}-a_2e^{-b_2}|\leq |a_1-a_2|+|a_1+a_2||b_1-b_2|, 
$$
for $a_1,a_2\in\mathbb R$ and $b_1, b_2\in\mathbb R_+$,  
for fixed $(s,x)\in \bar H_T$, 
$\alpha\in\A$, 
$\tau\in\ST(T-s)$  
and  
$\gamma=(\alpha,\eta,\xi)\in \A^{\varepsilon}$  we have 
$$
|w^{\gamma,\tau}(s,x)-w^{\alpha,\tau}(s,x)|
\leq N_0(I_1+I_2), 
$$
where
$$
I_1=E\int_0^{T-s}e^{-\lambda t}(1+t)
(|\varepsilon|+|x^{\gamma,s,x}_t-x^{\alpha,s,x}_t|)\,dt
$$
$$
\leq \varepsilon\int_0^Te^{-\lambda t}(t+1)\,dt
+3E\sup_{t\leq T-s}e^{-(\lambda-1)t}
|x^{\gamma,s,x}_t-x_t^{\alpha,s,x}|, 
$$
$$
I_2=Ee^{-\lambda\tau}
(|\varepsilon|+|x^{\gamma,s,x}_{\tau}-x_{\tau}^{\alpha,s,x}|)+
Ee^{-\lambda\tau}
\int_0^{\tau}(|\varepsilon|
+|x^{\gamma,s,x}_t-x_t^{\alpha,s,x}|)\,dt
$$
$$
\leq 2|\varepsilon|+2E\sup_{t\leq T-s}e^{-(\lambda-1)t}
|x^{\gamma,s,x}_t-x_t^{\alpha,s,x}|, 
$$
and $N_0$ is a constant depending only on $K$. 
By 
It\^o's formula we get 
$$
e^{-2(K^2+1)t}|x^{\gamma,s,x}_t-x_t^{\alpha,s,x}|^2
\leq N_0\int_0^te^{-2(K^2+1)r}\varepsilon^2\,dr+m_t
\leq N_1\varepsilon^2+m_t, 
$$
where $m$ is a local martingale and $N_1$ is a 
constant depending only on $K$. 
Hence 
$$
Ee^{-2(K^2+1)\rho}
|x^{\gamma,s,x}_{\rho}-x_{\rho}^{\alpha,s,x}|^2
\leq N_1\varepsilon^2
$$
for stopping times $\rho$, that by virtue 
of Lemma 3.2 from \cite{GK2} yields 
$$
E\sup_{t\leq T-s}e^{-(K^2+1)t}
|x^{\gamma,s,x}_{t}-x_{t}^{\alpha,s,x}|
\leq 3\sqrt{N_1}|\varepsilon|. 
$$
Consequently, \eqref{7.18.3.9} holds 
with a constant $N$ depending only on 
$K$ and $T$, and if $\lambda\geq K^2+2$ 
then $N$ is independent of $T$.  
\end{proof}

Now we shake the finite difference problem 
\eqref{1.18.02.09}-\eqref{2.18.02.09} when 
$m^{\alpha}=1$ for all $\alpha\in A$. 
We keep the notation of 
Section \ref{section-existence-of-soln-to-disc-prob} 
and Assumption \ref{assumption 1.23.12.08} in force.  
Moreover we make the following assumption. 

\begin{assumptions}                            \label{assumption 1.19.3.9}
For $\psi:=\sqrt{a^{\alpha}_k}, 
b^{\alpha}_k, f^{\alpha}, g,c^{\alpha}-\lambda$, 
for  
$\alpha\in A$ and $k=\pm1,\dots,\pm d_1$  
we have 
$$
|\psi(t,x)|\leq K, 
\quad
|\psi(t,x)-\psi(s,y)|\leq K(|x-y|+|s-t|^{1/2})
$$ 
for all $s,t\in\mathbb R$ and $x,y\in\mathbb R^d$. 
\end{assumptions}

Shaking the problem 
\begin{equation}                        \label{2.20.3.9}
\max[\sup_{\alpha\in A} 
\delta^T_{\tau}u+L_h^{\alpha} u+f^{\alpha}, g-u]=0 
\quad \text{on $H_T$}, 
\end{equation}
\begin{equation}                                  \label{3.20.3.9}
u(T,x)=g(T,x) \quad\text{for $x\in\mathbb R^d$}
\end{equation}
means that we consider the problem 
\begin{equation}                                  \label{4.20.3.9}
\max[\sup_{\gamma\in A^{\varepsilon}}
\delta^T_{\tau}u+L_h^{\gamma}u
+f^{\gamma}, g^{\varepsilon}-u]=0
\end{equation}
\begin{equation}                                   \label{5.20.3.9}
u(T,x)=g^{\varepsilon}(T,x)
\quad 
\text{for $x\in\mathbb R^d$}, 
\end{equation} 
where $g^{\varepsilon}$ is defined as in 
\eqref{10.19.3.9} and  for 
$\gamma=(\alpha,\xi,\eta,r)\in \bar A^{\varepsilon}$
$$
L_h^{\gamma}:=a_k^{\gamma}\Delta_{h,\ell_k}
+b_k^{\gamma}\delta_{h,\ell_k}-c^{\gamma},
$$
$c^{\gamma}$ and $f^{\gamma}$ are defined as in 
\eqref{13.19.3.9}, and 
$$
a_k^{\gamma}(t,x)=a_k^{\alpha}(t+\eta,x+\xi), 
\quad
b_k^{\gamma}(t,x)=b_k^{\alpha}(t+\eta,x+\xi),
\quad t\in\mathbb R, x\in\mathbb R^d,
$$
for $k=\pm1,\dots,\pm d_1$.

By virtue of Theorem \ref{theorem 1.20.3.9}, if 
Assumptions \ref{assumption 1.23.12.08} 
and \ref{assumption 1.19.3.9} hold then 
\eqref{2.20.3.9}-\eqref{3.20.3.9} and 
\eqref{4.20.3.9}-\eqref{5.20.3.9} have a unique 
bounded solution 
$w_{\tau,h}$ and $w^{\varepsilon}_{\tau,h}$, 
respectively.

\begin{lemma}                              \label{lemma 1.20.3.9}
Let Assumptions \ref{assumption 1.23.12.08} 
and \ref{assumption 1.19.3.9} hold. Then 
\begin{equation}                              \label{8.20.3.9}
|w^{\varepsilon}_{\tau,h}-w_{\tau,h}|\leq
N_0|\varepsilon|
\quad \text{on $\bar H_T$}, 
\end{equation}
with a constant $N_0$ depending only on $K$, $d_1$ and $T$. 
Assume, additionally, $\tau,h\leq K$. Then 
\begin{equation}                            \label{9.20.3.9}                              
|w^{\varepsilon}_{\tau,h}(s,x)-w_{\tau,h}^{\varepsilon}(t,y)|
\leq N_1(|x-y|+|s-t|^{1/2}+\sqrt{\tau})
\end{equation}
for all $s,t\in[0,T]$ and $x,y\in\mathbb R^d$, where 
$N_1$ is a constant depending only on $K$, $d_1$ and 
$T$. There is a constant $\lambda_0$ depending only on $K$ and 
$d_1$  
such that if $\lambda\geq\lambda_0$ then $N_0$ and $N_1$  
are independent of $T$. 
\end{lemma}

\begin{proof}
We get estimate \eqref{8.20.3.9} by an 
obvious application of Theorem \ref{theorem 3.17.3.9}. 
Estimate \eqref{9.20.3.9} follows immediately from 
Theorems \ref{theorem 4.17.3.9} and 
\ref{cor:hold:cont:in:t}. 
\end{proof}

Let $\rho\in C^{\infty}_0(\mathbb R^{d+1})$ 
be a fixed nonnegative function 
with support in $(-1,0)\times B_1$ and 
unit integral, where 
$B_1$ denotes the open ball of radius 1 
centered at the origin of $\mathbb R^d$.  
For $\varepsilon>0$ set 
$$
w^{\varepsilon(\varepsilon)}(t,x)
=\int_{\mathbb R^{d+1}}w^{\varepsilon}(s,y)
\rho((t-s)/\varepsilon^2, (x-y)/\varepsilon)\,ds\,dy 
$$ 
for $t\in[0,T]$ and $x\in\mathbb R^d$, 
where $w^{\varepsilon}(s,y):=w^{\varepsilon}(T,y)$ for 
$s\geq T$ and $y\in\mathbb R^d$. Define similarly 
$w^{\varepsilon(\varepsilon)}_{\tau,h}$ from  
$w^{\varepsilon}_{\tau,h}$. 

\begin{lemma}                                           \label{lemma 2.20.3.9} 
Let Assumption \ref{assumption 5.17.3.9} 
hold. Then there is a constant $N_0$ 
depending only on $K$ and $T$ such that 
\begin{equation}                                           \label{10.20.3.9}
|w^{\varepsilon(\varepsilon)}-w|
\leq
N_0\varepsilon\quad\text{on $\bar H_T$},  
\end{equation}
\begin{equation}                                           \label{11.18.3.9}
|w^{\varepsilon(\varepsilon)}(t,x)
-w^{\varepsilon(\varepsilon)}(s,y)|
\leq N_0(|x-y|+|t-s|^{1/2})
\end{equation}
for all $s,t\in[0,T]$ and $x,y\in\mathbb R^d$. 
For integers $n\geq1$ 
\begin{equation}                                            \label{1.21.3.9}
|D^n_tw^{\varepsilon(\varepsilon)}|+
|D^{2n}_xw^{\varepsilon(\varepsilon)}|
\leq N_1\varepsilon^{-2n+1}\quad
\text{on $\bar H_T$}, 
\end{equation}
where $N_1$ is a constant depending 
only on $n$, $K$, $d$ and $T$.  
There is a constant 
$\lambda_0$ such that if 
$\lambda\geq\lambda_0$ then $N_0$ and $N_1$ 
are independent 
of $T$. 
Moreover,
\begin{equation}                                             \label{5.21.3.9}
\max[D_tw^{\varepsilon(\varepsilon)}
+\sup_{\alpha\in A}(L^{\alpha}
w^{\varepsilon(\varepsilon)}+f^{\alpha}), 
g-w^{\varepsilon(\varepsilon)}]\leq 0
\quad \text{on $H_T$}. 
\end{equation}
\begin{proof}
Estimates 
\eqref{10.20.3.9}-\eqref{1.21.3.9} follow 
immediately from  
Lemma \ref{lemma 1.18.3.9}. 
To prove \eqref{1.21.3.9} we use 
\eqref{10.21.3.9} and define $w_n^{\varepsilon(\varepsilon)}$ 
from $w_n^{\varepsilon}$ as $w^{\varepsilon(\varepsilon)}$ 
is defined 
from $w^{\varepsilon}$. Notice that for $n\to\infty$
$$
D_tw^{\varepsilon(\varepsilon)}_n
\to w^{\varepsilon(\varepsilon)}, 
\quad D^{\beta}_xw^{\varepsilon(\varepsilon)}_n
\to
D^{\beta}_xw^{\varepsilon(\varepsilon)}_n
$$
for multi-indices $\beta$,  
by Lebesgue's theorem on dominated convergence. 
By Theorem 2.1 in \cite{krylov:rate:variable} 
for each integer $n\geq1$ we have 
$$
D_t w_n^{\varepsilon(\varepsilon)}
+L^{\alpha}w_n^{\varepsilon(\varepsilon)}
+f^{\alpha}+r(g^{\varepsilon}-w_n^{\varepsilon(\varepsilon)})
\leq0
\quad
\text{on $\bar H_T$}
$$
for all $\alpha\in A$ and $r\in[0,1]$. 
Letting here $n\to\infty$ and using that $g\leq g^{\varepsilon}$,  
we get 
$$
D_t w^{\varepsilon(\varepsilon)}
+L^{\alpha}w^{\varepsilon(\varepsilon)}
+f^{\alpha}+r(g-w^{\varepsilon(\varepsilon)})
\leq0
\quad
\text{on $\bar H_T$, for $\alpha\in A$, $r\geq0$}, 
$$
that is equivalent to \eqref{1.21.3.9}. 
\end{proof}
\end{lemma}

\begin{lemma}                                           \label{lemma 2.20.3.9} 
Let Assumptions \ref{assumption 1.23.12.08} 
and \ref{assumption 1.19.3.9} hold. 
Then, provided $T>2\varepsilon^2$, 
\begin{equation}                                         \label{12.21.3.9}
\max[\delta^T_{\tau}
w^{\varepsilon(\varepsilon)}_{\tau,h}
+\sup_{\alpha\in A}(L^{\alpha}_h
w^{\varepsilon(\varepsilon)}_{\tau,h}+f^{\alpha}), 
g-w^{\varepsilon(\varepsilon)}_{\tau,h}]\leq 0
\quad \text{on $H_{T-2\varepsilon^2}$}. 
\end{equation}
Assume, additionally, 
$\tau, h\leq K$. Then 
\begin{equation}                                         \label{4.21.3.9}
|w^{\varepsilon(\varepsilon)}_{\tau,h}-w_{\tau,h}|
\leq
N_0(|\varepsilon|+\sqrt{\tau})
\quad \text{on $\bar H_T$}, 
\end{equation}
\begin{equation}                               \label{12.20.3.9}   
|w^{\varepsilon(\varepsilon)}_{\tau,h}(t,x)
-w_{\tau,h}^{\varepsilon(\varepsilon)}(s,y)|
\leq N_0(|x-y|+|s-t|^{1/2}+\sqrt{\tau}), 
\end{equation}
for $t,s\in[0,T]$ and $x,y\in\mathbb R^d$, 
 where $N_0$ is a constant depending only on $K$,
$d_1$ and
$T$. Moreover, for $n\geq1$ there is a constant $N_1$ depending 
only on $n$, $K$, $d_1$, $d$ and $T$, such that 
\begin{equation}                                            \label{11.21.3.9}
|D^n_tw^{\varepsilon(\varepsilon)}_{\tau,h}|+
|D^{2n}_xw^{\varepsilon(\varepsilon)}_{\tau,h}|
\leq N_1\varepsilon^{-2n}(|\varepsilon|+\sqrt{\tau})
\quad
\text{on $\bar H_{T}$}.  
\end{equation}
There is a constant 
$\lambda_0$ depending on $K$ and $d_1$ such that 
if $\lambda\geq\lambda_0$ then $N_0$ and $N_1$ are independent 
of $T$. 
\end{lemma}
\begin{proof}
Estimates \eqref{4.21.3.9}-\eqref{11.21.3.9} follow immediately 
from Lemma \ref{lemma 1.20.3.9}. 
To prove \eqref{12.21.3.9} notice that from \eqref{4.20.3.9} we have 
for $\alpha\in A$ 
$$
(\delta^T_{\tau}+L^{\alpha}_h)
w^{\varepsilon(\varepsilon)}(t-\varepsilon^2s,x-\varepsilon y)
+f^{\alpha}(t,x)\leq 0, 
$$
$$
g(t,x)-w^{\varepsilon(\varepsilon)}(t-\varepsilon^2s,x-\varepsilon y)\leq0
$$
for $(t,x)\in \bar H_{T-2\varepsilon^2}$, $s\in[-1,0]$, 
$|y|\leq1$. Multiplying these inequalities by $\rho(s,y)$ and then integrating 
them against $ds\,dy$ we get \eqref{12.21.3.9}.  
\end{proof}

{\it Proof of Theorem \ref{theorem-the-main-result}:}
 
Let $\varepsilon = (\tau + h^2)^{1/4}$.  
Due to Theorems \ref{theorem 1.26.2.9} and 
\ref{cor:hold:cont:in:t} 
it suffices to consider the  
case $T>2\varepsilon^2$ and to prove 
\eqref{5.20.01.09} on $H_S$ with $S=T-2\varepsilon^2$. 
Notice that due to $\tau\leq 1$ we have 
$\tau<\varepsilon^2$. Hence for 
$u:=w^{\varepsilon(\varepsilon)},
w^{\varepsilon(\varepsilon)}_{\tau,h}$ we have 
$\delta^T_{\tau}u=\delta_{\tau}u$ on $H_S$, and  
by Taylor's formula and using 
\eqref{1.21.3.9} and \eqref{11.21.3.9}
$$
|\delta^T_{\tau}u-D_tu|+\sup_{\alpha\in A}|L^{\alpha}u-L^{\alpha}_hu|
$$
$$
\leq 
N_0(\tau\sup_{H_S}|D_t^2u|+h^2\sup_{H_S}|D^4_xu|
+h\sup_{H_S}|D^2_xu|)\leq N_1\varepsilon
$$   
on $H_S$. 
Notice also that 
$$
\sup_{H_T\setminus H_S}
(w_{\tau,h}-w^{\varepsilon(\varepsilon)})_+
$$
\begin{equation}                                              \label{15.21.3.9}
\leq \sup_{H_T\setminus H_S}(|w_{\tau,h}-g|+|g-w|+
|w-w^{\varepsilon(\varepsilon)}|)
\leq N_2\varepsilon, 
\end{equation}
$$
\sup_{\{S\}\times\mathbb R^d}
(w-w^{\varepsilon(\varepsilon)}_{\tau,h})_+
$$
\begin{equation}                                              \label{16.21.3.9}
\leq \sup_{\{S\}\times\mathbb R^d}(|w-g|+
|g-w_{\tau,h}|+
|w_{\tau,h}-w^{\varepsilon(\varepsilon)}_{\tau,h}|)
\leq N_2\varepsilon
\end{equation}
Thus by \eqref{5.21.3.9} 
 for $\alpha\in A$ 
\begin{equation}                               \label{1.22.3.9}                           
\delta_{\tau}^T\bar w^{\varepsilon(\varepsilon)}
+L_h^{\alpha}\bar w^{\varepsilon(\varepsilon)}+f^{\alpha}\leq 0 
\quad 
\text{and \,\,$g-\bar w^{\varepsilon(\varepsilon)}\leq0$ 
\,\,on $H_S$}
\end{equation}
\begin{equation}                              \label{2.22.3.9}
w_{\tau,h}\leq \bar w^{\varepsilon(\varepsilon)} 
\quad\text{on $\bar H_T\setminus H_S$}, 
\end{equation}
for 
$
\bar w^{\varepsilon(\varepsilon)}
:=w^{\varepsilon(\varepsilon)}+N_1(S-t)\varepsilon+N_2\varepsilon
$. 
If $\lambda>0$ then \eqref{1.22.3.9}  and \eqref{2.22.3.9} 
hold also for 
$\bar w^{\varepsilon(\varepsilon)}
:=w^{\varepsilon(\varepsilon)}+(N_1\lambda^{-1}+N_2)\varepsilon$.  
Similarly, by \eqref{12.21.3.9} for $\alpha\in A$ 
\begin{equation}                              \label{16.21.3.9}
D_t\bar w^{\varepsilon(\varepsilon)}_{\tau,h}
+L^{\alpha}\bar w^{\varepsilon(\varepsilon)}_{\tau,h}+f^{\alpha}
\leq0
\quad 
\text{and \,\,$g-\bar w^{\varepsilon(\varepsilon)}_{\tau,h}\leq0$ 
\,\,on $H_S$} 
\end{equation}
\begin{equation}                             \label{1.24.3.9}
w\leq \bar w^{\varepsilon(\varepsilon)}_{\tau,h}
\quad \text{on $\{S\}\times\mathbb R^d$}
\end{equation}
for   
$\bar w^{\varepsilon(\varepsilon)}_{\tau,h}
:=w^{\varepsilon(\varepsilon)}_{\tau,h}
+N_1(S-t)\varepsilon+N_2\varepsilon$ and also for 
$\bar w^{\varepsilon(\varepsilon)}_{\tau,h}
:=w^{\varepsilon(\varepsilon)}_{\tau,h}
+(N_1\lambda^{-1}+N_2)\varepsilon$ 
when $\lambda>0$. 
By Corollary \ref{corollary 1.25.2.9} from 
\eqref{1.22.3.9}-\eqref{2.22.3.9} we get 
$w_{\tau,h}\leq \bar w^{\varepsilon(\varepsilon)}$, and 
by Lemma \ref{lemma  1.11.02.09}
from \eqref{16.21.3.9}-\eqref{1.24.3.9} 
we have 
$w\leq \bar w^{\varepsilon(\varepsilon)}_{\tau,h}$ 
on $H_S$. Consequently, there is a constant $N$ such that 
$$
w_{\tau,h}\leq w+N\varepsilon, 
\quad
w\leq w_{\tau,h}+N\varepsilon \quad \text{on $H_S$}, 
$$
that obviously yields \eqref{5.20.01.09} on $H_S$. 
Inspecting the constants $N_0$, $N_1$ and $N_2$ we see that 
$N$ depends only on $K$, $d$, $d_1$ and $T$, and that there is a constant 
$\lambda_0$, depending only on $K$ and $d_1$ such that 
if $\lambda\geq\lambda_0$ then $N$ is independent of $T$.

\noindent
{\bf Acknowledgment.} 
The authors are grateful to Nicolai Krylov in Minnesota for
valuable information on the subject of this paper. 
They would like to thank the referee for noticing some mistakes 
and for useful suggestions. 

\chead{\textsc{References}}


\begin{thebibliography}{10}

\bibitem{barles:jakobsen:rate:hamilton}
Barles, G. and Jakobsen, E.~R. (2002).
\newblock On the convergence rate of approximation schemes for
  {H}amilton-{J}acobi-{B}ellman equations.
\newblock \emph{M2AN Math. Model. Numer. Anal.}, 36(1), 33--54.

\bibitem{barles:jakobsen:error:bounds}
Barles, G. and Jakobsen, E.~R. (2005).
\newblock Error bounds for monotone approximation schemes for
  {H}amilton-{J}acobi-{B}ellman equations.
\newblock \emph{SIAM J. Numer. Anal.}, 43(2), 540--558 (electronic).

\bibitem{biswas:jakobsen:karlsen:error}
Biswas, I.~H., Jakobsen, E.~R. and Karlsen, K.~H. (2006).
\newblock Error estimates for finite difference-quadrature schemes for a class
  of nonlocal {B}ellman equations with variable diffusion.
\newblock \emph{http://www.math.uio.no/eprint/pure\_math/2006/pure\_2006.html}.

\bibitem{dong:krylov:rate:domains}
Dong, H. and Krylov, N. (2007).
\newblock {The Rate of Convergence of Finite-Difference Approximations for
  Parabolic Bellman Equations with Lipschitz Coefficients in Cylindrical
  Domains}.
\newblock \emph{Applied Mathematics and Optimization}, 56(1), 37--66.
 

\bibitem{GK2} Gy\"ongy, I. and
Krylov, N. (2003).
On the rate of convergence of 
splitting-up 
approximations for SPDEs.  
In Progress in Probability, {\bf56}, 
Birkhauser Verlag, Basel, pp 301--321.

\bibitem{gyongy:siska:on:randomized}
Gy\"{o}ngy, I. and \v{S}i\v{s}ka, D. (2008).
\newblock On randomized stopping.
\newblock \emph{Bernoulli}, 14(2), 352--361.

\bibitem{jakobsen:rate:optimal:stopping}
Jakobsen, E.~R. (2003).
\newblock On the rate of convergence of approximation schemes for {B}ellman
  equations associated with optimal stopping time problems.
\newblock \emph{Math. Models Methods Appl. Sci.}, 13(5), 613--644.

\bibitem{jakobsen:karlsen:convergence:source:terms}
Jakobsen, E.~R. and Karlsen, K.~H. (2005).
\newblock Convergence rates for semi-discrete splitting approximations for
  degenerate parabolic equations with source terms.
\newblock \emph{BIT}, 45(1), 37--67.

\bibitem{jakobsen:karlsen:chioma:error}
Jakobsen, E.~R., Karlsen, K.~H. and La~Chioma, C. (2005).
\newblock Error estimates for approximate solutions to {B}ellman equations
  associated with controlled jump-diffusions.
\newblock \emph{http://www.math.uio.no/eprint/pure\_math/2005/pure\_2005.html}.

\bibitem{krylov:controlled}
Krylov, N.~V. (1980).
\newblock \emph{Controlled diffusion processes}, volume~14 of
  \emph{Applications of Mathematics}.
\newblock Springer-Verlag, New York.
\newblock Translated from the Russian by A. B. Aries.

\bibitem{krylov:rate:equations}
Krylov, N.~V. (1997).
\newblock On the rate of convergence of finite-difference approximations for
  {B}ellman's equations.
\newblock \emph{Algebra i Analiz}, 9(3), 245--256.

\bibitem{krylov:approximating:value}
Krylov, N.~V. (1999).
\newblock Approximating value functions for controlled degenerate diffusion
  processes by using piece-wise constant policies.
\newblock \emph{Electronic Journal of Probability}, 4(2), 1--19.

\bibitem{krylov:rate:variable}
Krylov, N.~V. (2000).
\newblock On the rate of convergence of finite-difference approximations for
  {B}ellman's equations with variable coefficients.
\newblock \emph{Probab. Theory Related Fields}, 117(1), 1--16.

\bibitem{krylov:rate:lipschitz}
Krylov, N.~V. (2004).
\newblock On the rate of convergence of finite-difference approximations for
  {B}ellman equations with {L}ipschitz coefficients.
\newblock \emph{arXiv:math}, 1(1), 1--33.

\bibitem{krylov:rate:lipschitz:published}
Krylov, N.~V. (2005).
\newblock The rate of convergence of finite-difference approximations for
  {B}ellman equations with {L}ipschitz coefficients.
\newblock \emph{Appl. Math. Optim.}, 52(3), 365--399.

\bibitem{krylov:factorizations}
Krylov, N.~V. (2008).
\newblock On factorizations of smooth nonnegative matrix-values functions and
  on smooth functions with values in polyhedra.
\newblock \emph{Appl. Math. Optim.}, 58(3), 373--392.


\bibitem{krylov:apriori}
Krylov, N.~V. (2007).
\newblock A priori estimates of smoothness of solutions to difference {B}ellman
  equations with linear and quasi-linear operators.
\newblock \emph{Math. Comp.}, 76(258), 669--698.

\bibitem{kushner:dupuis:numerical}
Kushner, H.~J. and Dupuis, P. (2001).
\newblock \emph{Numerical methods for stochastic control problems in continuous
  time}, volume~24 of \emph{Applications of Mathematics (New York)}.
\newblock Springer-Verlag, New York, second edition.
\newblock Stochastic Modelling and Applied Probability.

\bibitem{menaldi:some:estimates}
Menaldi, J.-L. (1989).
\newblock Some estimates for finite difference approximations.
\newblock \emph{SIAM J. Control Optim.}, 27(3), 579--607.

\bibitem{shiryaev:statistical}
Shiryaev, A.~N. (1976).
\newblock \emph{Statisticheskii posledovatelnyi analiz. {O}ptimalnye pravila
  ostanovki}.
\newblock Izdat. ``Nauka'', Moscow.
\newblock Second edition, revised.

\end{thebibliography}
\end{document}